\documentclass[12pt]{amsart}
\usepackage[T1]{fontenc}
\usepackage{amsmath}
\usepackage{amsfonts}
\usepackage{amssymb}
\usepackage{xcolor}
\usepackage{bbding}
\usepackage{makecell}
\usepackage{stmaryrd}
\usepackage{txfonts}
\usepackage{graphicx}
\usepackage{diagbox}
\usepackage{epsfig}
\usepackage{xypic}
\usepackage{tikz}
\usepackage{longtable}
\usepackage{caption}
\usepackage[all]{xy}
\usepackage{pgflibraryarrows}
\usepackage{pgflibrarysnakes}
\usepackage[shortlabels]{enumitem}
\usepackage{ifpdf}
\ifpdf
  \usepackage[colorlinks,final,hyperindex]{hyperref}
\else
  \usepackage[colorlinks,final,hyperindex,driverfallback=hypertex]{hyperref}
\fi

\usepackage{xspace}

\usepackage{listings}
\newtheorem{theorem}{Theorem}[section]
\newtheorem{lemma}[theorem]{Lemma}
\newtheorem{coro}[theorem]{Corollary}

\newtheorem{conjecture}[theorem]{Conjecture}

\theoremstyle{definition}

\newtheorem{remark}[theorem]{Remark}
\newtheorem{exmp}[theorem]{Example}

\newcommand{\nc}{\newcommand}
\newcommand{\delete}[1]{}

\nc{\tred}[1]{\textcolor{red}{#1}}
\nc{\tblue}[1]{\textcolor{blue}{#1}} \nc{\tgreen}[1]{\textcolor{green}{#1}} \nc{\tpurple}[1]{\textcolor{purple}{#1}} \nc{\btred}[1]{\textcolor{red}{\bf #1}} \nc{\btblue}[1]{\textcolor{blue}{\bf #1}} \nc{\btgreen}[1]{\textcolor{green}{\bf #1}} \nc{\btpurple}[1]{\textcolor{purple}{\bf #1}}

\newcommand{\efootnote}[1]{}
\nc{\mlabel}[1]{\label{#1}}  % Use this to suppress names
\nc{\mcite}[2][]{\cite[#1]{#2}}  % Use this to suppress names
\nc{\mref}[1]{\ref{#1}}  % Use this to suppress names
\nc{\mbibitem}[1]{\bibitem{#1}} % Use this to show number
\delete{
\nc{\mlabel}[1]{\label{#1}  % Use the next two lines to show names
{\hfill \hspace{1cm}{\bf{{\ }\hfill(#1)}}}}
\nc{\mcite}[1]{\cite{#1}}  % Use this lines to show names
\nc{\mref}[1]{\ref{#1}{{\bf{{\ }(#1)}}}}  % Use this lines to show names
\nc{\mbibitem}[1]{\bibitem[\bf #1]{#1}} % Use this to show name
}

\renewcommand\geq{\geqslant}
\renewcommand\leq{\leqslant}

\renewcommand\bar[1]{\overline{#1}}
\nc{\nz}{\varepsilon}
\nc{\Id}{\mathrm{Id}}
\nc{\map}[2]{{#2}^{#1}}
\nc{\gp}{B}

\nc{\Irr}{\mathrm{Irr}}
\nc{\vx}{\sigma} \nc{\vy}{\tau} \nc{\dvx}{\sigma^{(1)}} \nc{\dvy}{\tau^{(1)}} \nc{\done}{\vep} \nc{\mcitep}[1]{\mcite{#1}} \nc{\wt}{\mathrm{wt}} \nc{\bre}[1]{|#1|} \nc{\mapmonoid}{\frakM} \nc{\disjoint}{\frakM'}
\nc{\ncpoly}[1]{\langle #1\rangle}  %for noncommutative polynomials
\nc{\mapm}[1]{\lfloor\!|{#1}|\!\rfloor}
\nc{\diff}[1]{{}^\NC\{ #1 \}} \nc{\disj}[1]{\{{#1}\}'} \nc{\mdisj}[1]{\frakM'(#1)} \nc{\brho}{\bar{\rho}} \nc{\om}{\bar{\frakm}} \nc{\frakn}{\mathfrak n} \nc{\ddeg}[1]{^{(#1)}} \nc{\opset}{X} \nc{\genset}{{Z}} \nc{\NC}{\mathrm{{NC}}} \nc{\leaf}{\mathrm{leaf}} \nc{\twig}{\mathrm{twig}} \nc{\fe}{\mathrm{fl}} \nc{\munderline}[1]{#1} \nc{\bo}{o} \nc{\dep}{\mathrm{depth}} \nc{\ofe}{\mathrm{ofl}} \nc{\dfe}{\mathrm{dfe}} \nc{\fex}{\mathrm{fex}} \nc{\dl}{\mathrm{dlex}} \nc{\db}{\mathrm{db}} \nc{\lex}{\mathrm{lex}} \nc{\clex}{\mathrm{clex}} \nc{\dgp}{\mathrm{dgp}} \nc{\dgx}{\mathrm{dgx}} \nc{\br}{\mathrm{br}} \nc{\obd}{\mathrm{odb}} \nc{\ob}{\mathrm{ob}}
\nc{\pie}{\mathrm{PIE}}
\nc{\rbo}{\mathrm{RBO}}
\nc{\supp}{\mathcal{S}}
\nc{\nul}{\mathcal{Z}}
\nc{\bin}[2]{ (_{\stackrel{\scs{#1}}{\scs{#2}}})}  %binomial coeff
\nc{\binc}[2]{ \left (\!\! \begin{array}{c} \scs{#1}\\
    \scs{#2} \end{array}\!\! \right )}  %binomial coeff
\nc{\bincc}[2]{  \left ( {\scs{#1} \atop
    \vspace{-1cm}\scs{#2}} \right )}  %binomial coeff
\nc{\bs}{\bar{S}} \nc{\cosum}{\sqsubset} \nc{\la}{\longrightarrow} \nc{\rar}{\rightarrow} \nc{\dar}{\downarrow} \nc{\dprod}{**} \nc{\dap}[1]{\downarrow \rlap{$\scriptstyle{#1}$}} \nc{\md}[1]{\bar{#1}} \nc{\uap}[1]{\uparrow \rlap{$\scriptstyle{#1}$}} \nc{\defeq}{\stackrel{\rm def}{=}} \nc{\disp}[1]{\displaystyle{#1}} \nc{\dotcup}{\ \displaystyle{\bigcup^\bullet}\ } \nc{\gzeta}{\bar{\zeta}} \nc{\hcm}{\ \hat{,}\ } \nc{\hts}{\hat{\otimes}} \nc{\barot}{{\otimes}} \nc{\free}[1]{\bar{#1}} \nc{\uni}[1]{\tilde{#1}} \nc{\hcirc}{\hat{\circ}} \nc{\leng}{\ell} \nc{\lleft}{[} \nc{\lright}{]} \nc{\lc}{\lfloor} \nc{\rc}{\rfloor}
\nc{\lb}{[} %left bracket
\nc{\rb}{]} %right bracket
\nc{\curlyl}{\left \{ \begin{array}{c} {} \\ {} \end{array}
    \right.  \!\!\!\!\!\!\!}
\nc{\curlyr}{ \!\!\!\!\!\!\!
    \left. \begin{array}{c} {} \\ {} \end{array}
    \right \} }
\nc{\longmid}{\left | \begin{array}{c} {} \\ {} \end{array}
    \right. \!\!\!\!\!\!\!}
\nc{\onetree}{\bullet} \nc{\ora}[1]{\stackrel{#1}{\rar}}
\nc{\ola}[1]{\stackrel{#1}{\la}}%${\Bbb Z}$
\nc{\ot}{\otimes} \nc{\mot}{{{\boxtimes\,}}} \nc{\otm}{\overline{\boxtimes}} \nc{\sprod}{\bullet} \nc{\scs}[1]{\scriptstyle{#1}} \nc{\mrm}[1]{{\rm #1}} \nc{\msum}{\sum\limits}
\nc{\margin}[1]{\marginpar{\rm #1}}   %{\rm #1}}
\nc{\dirlim}{\displaystyle{\lim_{\longrightarrow}}\,} \nc{\invlim}{\displaystyle{\lim_{\longleftarrow}}\,} \nc{\mvp}{\vspace{0.3cm}} \nc{\tk}{^{(k)}} \nc{\tp}{^\prime} \nc{\ttp}{^{\prime\prime}} \nc{\svp}{\vspace{2cm}} \nc{\vp}{\vspace{8cm}} \nc{\proofbegin}{\noindent{\bf Proof: }}
\nc{\proofend}{$\blacksquare$ \vspace{0.3cm}}
\nc{\modg}[1]{\!<\!\!{#1}\!\!>}
\nc{\intg}[1]{F_C(#1)} \nc{\lmodg}{\!<\!\!} \nc{\rmodg}{\!\!>\!} \nc{\cpi}{\widehat{\Pi}}
\nc{\sha}{{\,\makebox[0.6em]{\cyr X}\,}}  %used to be \cyr
\nc{\shap}{{\mbox{\cyrs X}}} %sha as product
\nc{\shpr}{\diamond}    %Shuffle product
\nc{\shp}{\ast} \nc{\shplus}{\shpr^+}
\nc{\shprc}{\shpr_c}    %Cartier's product
\nc{\msh}{\ast} \nc{\zprod}{m_0} \nc{\oprod}{m_1} \nc{\vep}{\varepsilon} \nc{\labs}{\mid\!} \nc{\rabs}{\!\mid}
\nc{\astarrow}{\overset{\raisebox{-3pt}{$\ast$}}{\rightarrow}}
\nc{\Sym}{{\mrm{Sym}}}
\nc{\Nsym}{{\mrm{NSym}}}
\nc{\Qsym}{{\mrm{QSym}}}
\nc{\SSym}{{\mrm{\mathfrak{S}Sym}}}
\nc{\SSRCT}{{\mrm{SSRCT}}}
\nc{\SSRCTs}{{\mrm{SSRCTs}}}
\nc{\SRCT}{{\mrm{SRCT}}}
\nc{\SRCTs}{{\mrm{SRCTs}}}
\nc{\SSRT}{{\mrm{SSRT}}}
\nc{\SSRTs}{{\mrm{SSRTs}}}
\nc{\syms}{{symmetric functions\xspace}}
\nc{\qsyms}{{quasi-symmetric functions\xspace}}
\nc{\nsymg}{{\mathrm{NSym}_\gp}}
\nc{\HSym}{{\mrm{\mathfrak{H}Sym}}}
\nc{\parr}{{\mrm {Par}}}
\nc{\Set}{{\mrm {Set}}}
\nc{\Comp}{\mrm {Comp}}
\nc{\comp}{\mrm {comp}}
\nc{\Des}{{\mrm {Des}}}
\nc{\pc}{\mrm {pc}}
\nc{\Sol}{{\mrm {Sol}}}
\nc{\Span}{\mrm {span}}
\nc{\Sh}{{\mrm {Sh}}}
\nc{\st}{{\rm{st}}}
\nc{\dth}{d} \nc{\mmbox}[1]{\mbox{\ #1\ }} \nc{\fp}{\mrm{FP}} \nc{\rchar}{\mrm{char}} \nc{\Fil}{\mrm{Fil}} \nc{\Mor}{Mor\xspace} \nc{\gmzvs}{gMZV\xspace} \nc{\gmzv}{gMZV\xspace} \nc{\mzv}{MZV\xspace} \nc{\mzvs}{MZVs\xspace} \nc{\Hom}{\mrm{Hom}} \nc{\id}{\mrm{id}} \nc{\im}{\mrm{im}} \nc{\incl}{\mrm{incl}}  \nc{\mchar}{\rm char}

\nc{\Alg}{\mathbf{Alg}} \nc{\Bax}{\mathbf{Bax}} \nc{\bff}{\mathbf f} \nc{\bfk}{{\bf k}} \nc{\bfone}{{\bf 1}} \nc{\bfx}{\mathbf x} \nc{\bfy}{\mathbf y}
\nc{\base}[1]{\bfone^{\otimes ({#1}+1)}} %{{a_{#1}}}
\nc{\Cat}{\mathbf{Cat}} \delete{}
\nc{\detail}{\marginpar{\bf More detail}
    \noindent{\bf Need more detail!}
    \svp}
\nc{\Int}{\mathbf{Int}} \nc{\Mon}{\mathbf{Mon}}
\nc{\rbtm}{{shuffle }} \nc{\rbto}{{Rota-Baxter }} \nc{\remarks}{\noindent{\bf Remarks: }} \nc{\Rings}{\mathbf{Rings}} \nc{\Sets}{\mathbf{Sets}}
\nc{\balpha}{\mathbf{\alpha}}

\nc{\BA}{{\mathbb A}} \nc{\CC}{{\mathbb C}} \nc{\DD}{{\mathbb D}} \nc{\EE}{{\mathbb E}} \nc{\FF}{{\mathbb F}} \nc{\GG}{{\mathbb G}} \nc{\HH}{{\mathbb H}} \nc{\LL}{{\mathbb L}} \nc{\NN}{{\mathbb N}} \nc{\KK}{{\mathbb K}} \nc{\PP}{{\mathbb P}} \nc{\QQ}{{\mathbb Q}} \nc{\RR}{{\mathbb R}} \nc{\TT}{{\mathbb T}} \nc{\VV}{{\mathbb V}} \nc{\ZZ}{{\mathbb Z}}

\nc{\cala}{{\mathcal A}} \nc{\calc}{{\mathcal C}} \nc{\cald}{{\mathcal D}} \nc{\cale}{{\mathcal E}} \nc{\calf}{{\mathcal F}} \nc{\calg}{{\mathcal G}} \nc{\calh}{{\mathcal H}} \nc{\cali}{{\mathcal I}} \nc{\call}{{\mathcal L}} \nc{\calm}{{\mathcal M}} \nc{\caln}{{\mathcal N}} \nc{\calo}{{\mathcal O}} \nc{\calp}{{\mathcal P}} \nc{\calr}{{\mathcal R}} \nc{\cals}{{\mathcal S}} \nc{\calt}{{\mathcal T}} \nc{\calw}{{\mathcal W}} \nc{\calk}{{\mathcal K}} \nc{\calx}{{\mathcal X}}
\nc{\calz}{{\mathcal Z}}
 \nc{\CA}{\mathcal{A}}

\nc{\fraka}{{\mathfrak a}} \nc{\frakA}{{\mathfrak A}} \nc{\frakb}{{\mathfrak b}} \nc{\frakB}{{\mathfrak B}}
\nc{\frakc}{{\mathfrak c}}  \nc{\frakD}{{\mathfrak D}}
\nc{\frakH}{{\mathfrak H}}
\nc{\frakh}{{\mathfrak h}} \nc{\frakM}{{\mathfrak M}}
\nc{\frakO}{{\mathfrak O}}
\nc{\frakE}{{\mathfrak E}}
\nc{\bfrakM}{\overline{\frakM}} \nc{\frakm}{{\mathfrak m}} \nc{\frakP}{{\mathfrak P}} \nc{\frakN}{{\mathfrak N}} \nc{\frakp}{{\mathfrak p}} \nc{\frakS}{{\mathfrak S}}
\nc{\frakk}{{\mathfrak k}}
\nc{\frakx}{{\mathfrak x}}
\nc{\frakl}{{\mathfrak l}} \nc{\ox}{\bar{\frakx}} \nc{\frakX}{{\mathfrak X}} \nc{\fraky}{{\mathfrak y}} \nc\dop{\delta}
\nc{\Reduce}{{\rm Red}}

\nc{\name}[1]{{\bf #1}}

\font\cyr=wncyr10 \font\cyrs=wncyr7
\nc{\redt}[1]{\textcolor{red}{#1}}
\nc{\yu}[1]{\textcolor{red}{\tt Yu:#1}}
\nc{\gu}[1]{\textcolor{blue}{\tt Gu:#1}}
\begin{document}
\title[Cancellation-free antipode formula for the Malvenuto-Reutenauer Hopf Algebra]{On the cancellation-free antipode formula for the Malvenuto-Reutenauer Hopf Algebra}

\author{Da Xu}
\address{School of Mathematics and Statistics, Southwest University, Chongqing 400715, China}
\email{x15903878019@email.swu.edu.cn}

\author{Houyi Yu$^{\ast}$}  \thanks{*Corresponding author}
\address{School of Mathematics and Statistics, Southwest University, Chongqing 400715, China}
\email{yuhouyi@swu.edu.cn}
%========================================================================
\hyphenpenalty=8000
%\date{\today}

\begin{abstract}
For the Malvenuto-Reutenauer Hopf algebra of permutations, we provide a cancellation-free antipode formula for any permutation of the form %in the one-line notation
$ab1\cdots(b-1)(b+1)\cdots(a-1)(a+1)\cdots n$, which
starts with the decreasing sequence $ab$ and ends with the increasing sequence $1\cdots(b-1)(b+1)\cdots(a-1)(a+1)\cdots n$, where $1\leq b<a\leq n$. As a consequence, we confirm two conjectures posed by Carolina Benedetti and Bruce E. Sagan.
\end{abstract}

\keywords{Malvenuto-Reutenauer Hopf algebra, permutation, antipode, cancellation-free}

\maketitle

\hyphenpenalty=8000 \setcounter{section}{0}

%========================================================================

\allowdisplaybreaks
%====================ÒÔÏ¿ªÊ¼ÕýÎÄ==================
\section{Introduction}\label{sec:int}
\delete{
In a highly influential paper \cite{JR79}, Join and Rota observed that Hopf algebras provide a valuable formal framework for the study of assembly and
disassembly operations of combinatorial objects.
Since then numerous combinatorial objects have been endowed with connected graded Hopf algebra structures which, while capturing precise algebraic properties of combinatorial operations, provided intuitive explanations of abstract concepts. Central among these combinatorial Hopf algebras are the Hopf algebras of symmetric functions \cite{Ge77}, of quasi-symmetric functions \cite{Ge84} and of noncommutative symmetric functions \cite{GKLLRT95}, and the Malvenuto-Reutenauer Hopf algebra of permutations \cite{MR95}, and so on.
}

The antipode of a combinatorial Hopf algebra can be computed by the celebrated Takeuchi's formula \cite{MT71}.
However, this antipode formula usually contains a great number of cancellations in its alternating sum.
Recently, significant attention has been dedicated to developing cancellation-free antipode formulas for various combinatorial Hopf algebras, including the Hopf algebra of graphs, matroids, simplicial complexes, set operads and posets, among others \cite{BB19,BC17,BBM19,BHM16,BEJM18,EH19,HM12,ML14,Pa16}.
Notably, Benedetti and Sagan \cite{BS17} employed the approach of sign-reversing involutions to establish cancellation-free antipode formulas
for nine different combinatorial Hopf algebras.

The Malvenuto-Reutenauer Hopf algebra $\SSym$ of permutations, introduced by Malvenuto and Reutenauer \cite{Mal94,MR95},  is a connected graded Hopf algebra which has a linear basis consisting of permutations in all symmetric groups.
Aguiar and Sottile \cite{AS05} gave, for the first time, an explicit antipode formula for $\SSym$, but it involves massive cancellations. Then Aguiar and Mahajan \cite{AM10} provided a cancellation-free expression for the antipode in terms of Hopf monoids.
In the seminal paper \cite{BS17}, Benedetti and Sagan recovered the antipode formulas for certain permutations by using the technique of sign-reversing involutions and appealing only to the Hopf algebra structure.
In the same paper the authors predicted two cancellation-free antipode formulas for some permutations whose image under the Robinson-Schensted map is a column superstandard Young tableau of hook shape (see Conjectures 9.9 and 9.10 of \cite{BS17} for details).

The purpose of this work is to prove the two conjectures posed by Benedetti and Sagan. More precisely, we give a cancellation-free antipode formula for any permutation, in the one-line notation, of the form $ab1\cdots(b-1)(b+1)\cdots(a-1)(a+1)\cdots n$, which
starts with the decreasing sequence $ab$ and ends with the increasing sequence $1\cdots(b-1)(b+1)\cdots(a-1)(a+1)\cdots n$, where $1\leq b<a\leq n$.
The two conjectures then follow directly as consequences of our result.

The organization of the paper is as follows. In the next section, after recalling some definitions related to the Malvenuto-Reutenauer Hopf algebra,
we state the main result, Theorem \ref{maintheorem}. Then, by using mathematical induction, we prove our main result in Section \ref{Proofofthemainythm}, where the two conjectures of Benedetti and Sagan are proved as consequences of Theorem \ref{maintheorem}.

\section{Preliminaries}
We begin by making clear some notation and definitions concerning the Malvenuto-Reutenauer Hopf algebra, and referring the reader to \cite{GR14,MR95} for more details.
For an introduction to Hopf algebras, we refer the reader to \cite{Sw69}.

Throughout, let $\mathbb{N}$ and $\mathbb{P}$ denote the set of nonnegative integers and positive integers, respectively.
Given $n\in\mathbb{P}$, we let $[n]$ denote the set $\{1, 2, \ldots, n\}$.
Also, $\bfk$ will usually denote a field of characteristic zero, although it may also be a commutative ring.

For an alphabet $A$, let $A^*$ be the set of words on $A$, including the empty word.
Let $\bfk\langle A\rangle$ be the noncommutative algebra over $A$.
The \emph{shuffle product} $\sha$ on $\bfk\langle A\rangle$ is defined by making $1\in \bfk\langle A\rangle$ the identity element for each product, and requiring that
it satisfies the following recursive relation
\begin{align}\label{eq:aushabv}
au\sha bv=a(u\sha bv)+b(au\sha v)
\end{align}
for all $a,b\in A$  and $u,v\in A^*$.
The following equivalent formulation of Eq.\,\eqref{eq:aushabv} is sometimes convenient:
\begin{align}\label{eq:uashavb}
ua\sha vb=(u\sha vb)a+(ua\sha v)b,
\end{align}
where $a,b\in A$  and $u,v\in A^*$.
Here and throughout the paper, the concatenation product always takes precedence over the shuffle product. Hence, in Eq.\eqref{eq:aushabv},
$(au)\sha(bv)$ is simply written as $au\sha bv$ by omitting all of brackets.
For example, for any $a,b,c,d\in A$, we have
$$ab\sha cd=abcd+acbd+acdb+cabd+cadb+cdab.$$
A word appearing in the shuffle product of $u\sha v$ is called a \emph{shuffle} of $u$ and $v$.
It will be convenient to identify a shuffle set with the sum of its elements.

Let $\mathfrak{S}_n$ be the group of permutations of $[n]$ and $\mathfrak{S} = \bigcup_{n\geq0}\mathfrak{S}_n$, where $\mathfrak{S}_0$ consists of the empty permutation $\emptyset$.
The \emph{Malvenuto-Reutenauer Hopf algebra of permutations} \cite{MR95} is a connected graded Hopf algebra whose underlying space is the graded $\bfk$-vector space $\SSym:=\bigoplus_{n\geq0}\bfk\mathfrak{S}_n$
with basis $\mathfrak{S}$. The counit is the projection $\varepsilon:\SSym\rightarrow \bfk\mathfrak{S}_0$.
To describe its Hopf structure, we need more notions.

A permutation $\sigma\in \mathfrak{S}_n$ is usually written in the one-line notation $\sigma_1\sigma_2\cdots\sigma_n$, where $\sigma_i=\sigma(i)$ for all $i\in [n]$. Let $w=w_1w_2\cdots w_n$ be a sequence of positive integers. If $m$ is a positive integer, then let $w+m$ be the word obtained by replacing in $w$ each $i$ by $i+m$.
The \emph{standardization} of $w$, denoted by $\st(w)$, is the unique permutation $\sigma\in \mathfrak{S}_n$ such that
$$
\sigma_i<\sigma_j \Leftrightarrow w_i\leq w_j \quad \text{whenever}\quad 1\leq i<j\leq n.
$$
For example, if $w=113$, then $w+3=446$ and $\st(w)=123$.
In this paper, because we concern mainly with elements of $\SSym$, we only take standardization of words with no repeated characters.
Now if $\pi\in \mathfrak{S}_m$ and $\sigma\in\mathfrak{S}_n$, then the product $\overline{\sha}$ and coproduct $\Delta$ in $\SSym$ are defined by
\begin{align*}
    \pi\bar{\sha} \sigma=\sum_{\tau\,\in\,\pi{\sha}(\sigma+m)}\tau\quad\text{and}\quad \Delta(\pi)=\sum_{i=0}^m\st(\pi_1\cdots\pi_i)\otimes \st(\pi_{i+1}\cdots\pi_{m}),
\end{align*}
respectively. The product is often called the \emph{shifted shuffle product}.
To illustrate,
\begin{align*}
  21\bar{\sha}12 =& 2134 + 2314 + 2341 + 3214 + 3241 + 3421\\
    \intertext{and}
  \Delta(2314) =& \emptyset \otimes 2314 + 1 \otimes 213 + 12 \otimes 12 + 231 \otimes 1 + 2314 \otimes \emptyset.
\end{align*}
The antipode $S$ for $\SSym$ can be computed recursively by
\begin{align}\label{eq:antipodeformula1}
    S(\sigma)=-\sum_{i=0}^{n-1}S(\st(\sigma_1\cdots \sigma_i))\bar{\sha}\st(\sigma_{i+1}\cdots \sigma_{n}),
\end{align}
where $\sigma\in\mathfrak{S}_n$ with $n\geq0$.
For notational convenience, for any $\sigma\in\mathfrak{S}_n$ and $j\in [n]$, we write $S(\sigma)_j$ for the sum of the terms of $S(\sigma)$ ending in $j$ with signs,
and $S(\sigma)_j^*$ for the sum obtained from $S(\sigma)_j$ by omitting the common suffix $j$.
The notions of $M_j$ and $M_j^*$ are similarly defined for any nonempty shuffle set $M$.
Thus, we have
\begin{align*}
	S(\sigma)=\sum_{j=1}^nS(\sigma)_j=\sum_{j=1}^nS(\sigma)_j^*j.
\end{align*}
For example, by Eq.\,\eqref{eq:antipodeformula1}, we have $S(213)=231-132-312=231-(13+31)2$. Hence, $S(213)_1=231$, $S(213)_2=-132-312$, $S(213)_3=0$, $S(213)_1^*=23$, $S(213)_2^*=-13-31$,
and $S(213)_3^*=0$.

For integers $k$ and $l$, let
\begin{align*}%\label{eq:etakldeltalk}
    \eta_{k,l}=\begin{cases}
    k(k+1)\cdots l&\text{if}\ k\leq l,\\
    \emptyset&\text{if}\ k=l+1,\\
    0&\text{if}\ k\geq l+2,
    \end{cases}
    \qquad\text{and}\qquad \delta_{l,k}=\begin{cases}
    l(l-1)\cdots k&\text{if}\ k\leq l,\\
    \emptyset&\text{if}\ k=l+1,\\
    0&\text{if}\ k\geq l+2.
    \end{cases}
\end{align*}
Since we only deal with results for $\SSym$, any words containing $0$ will be considered to be the zero element of $\SSym$.
 For example, if $n=4$, then $13\sha42\sha\delta_{n,5}=13\sha42\sha\emptyset=13\sha42$ and $5412\sha\delta_{n,6}=5412\sha 0=0$. If $a=2$ and $n\geq4$, then $a(a-1)\eta_{1,a-2}\sha \delta_{n,a+2}=21\emptyset\sha \delta_{n,4}=21\sha \delta_{n,4}$.

More generally, for a set $A$ of positive integers, let $\eta_A$ and $\delta_A$ denote the increasing and decreasing words whose
elements are in $A$, respectively. Given $A\subseteq[n]$, we let $\overline{A}=[n]-A$ be the complement of $A$ in $[n]$, and write $\sigma_{A,n}=\delta_{A}\eta_{\bar{A}}$.
When no confusion will result we will simply write $\sigma_{A}=\sigma_{A,n}.$
%$$\eta_k=\eta_{1,k},\quad \delta_k=\delta_{n,k}\quad \text{and}\quad \sigma_{A}=\sigma_{A,n}.$$
For example, for $n\geq 5$ and $A=\{2,5\}$, we have $\delta_{A}=52$ and $\eta_{\bar{A}}=13467\cdots n$, so that $\sigma_{A}=5213467\cdots n$.
More generally, for $A=\{a,b\}$ with $1\leq b < a \leq n$, we have
\begin{align*}
	\sigma_{A} =ab\eta_{1,b-1}\eta_{b+1,a-1}\eta_{a+1,n}.
\end{align*}

The following two conjectures were made in \cite[Conjectures 9.9 and 9.10]{BS17}. Here and throughout the paper, we use the notation $\sum_{i=j}^ka_i=0$ if $k<j$.

\begin{conjecture}[See \cite{BS17}, Conjecture 9.9]\label{conj1}
Let $A=\{a\}$ with $1 < a \leq n$. We have
\begin{align*}
 S(\sigma_{A}) =& (-1)^{n-1}\left(2 \sha \delta_{n,4}\right)31 + (-1)^{n+a}\left((a-1)\eta_{1,a-2} \sha \delta_{n,a+1}\right)a \\
  &+ \sum_{j=2}^{a-1}(-1)^{n+j}\left[\left((j-1)\eta_{1,j-2} \sha \delta_{n,j+1}\right)j + \left((j+1)\eta_{1,j-1} \sha \delta_{n,j+2}\right)j\right].
\end{align*}
\end{conjecture}

\begin{conjecture}[See \cite{BS17}, Conjecture 9.10]\label{conj2}
Let $A=\{a, 2\}$ with $2 < a \leq n$. We have
\begin{align*}
 S(\sigma_{A}) =& (-1)^{n}\left[\left(32 \sha \delta_{n,5}\right)41 + \left(12 \sha \delta_{n,4}\right)3\right] + (-1)^{n-1}\left[\left(1 \sha (3 \sha \delta_{n,5}\right)4)2\right]\\
 &+ \sum_{j=3}^{a-1} (-1)^{n+j}\left[\left((j+1)21\eta_{3,j-1} \sha \delta_{n,j+2}\right)j - \left(j21\eta_{3,j-1} \sha \delta_{n,j+2}\right)(j+1)\right].
\end{align*}
\end{conjecture}

Let $A=\{a,b\}$ with $1\leq b < a \leq n$.
Note that $\sigma_{\{a\}}=\sigma_{\{1,a\}}$. We will establish the following cancellation-free antipode formula for permutations of the type $\sigma_{A}$, and thereby verify
the two conjectures as corollaries.
\begin{theorem}\label{maintheorem}
Let $A=\{a,b\}$ and $\sigma_{A}\in \mathfrak{S}_n$, where $a,b$ and $n$ are positive integers satisfying $1 \leq b < a \leq n$.
\begin{enumerate}
  \item\label{lem:item:j=1b=1}
If $j=1$ and $b=1$, then $S(\sigma_{A})_j^*=(-1)^{n+1}(2\sha \delta_{n,4})3$.

\item\label{lem:item:j=1b=2}
If $j=1$ and $b=2$, then $S(\sigma_{A})_j^*=(-1)^{n} (32\sha \delta_{n,5})4$.

\item\label{lem:item:j=1bgeq3}
If $j=1$ and  $b\geq3$, then
$$S(\sigma_{A})_j^*=(-1)^{n+b+1}\left[(b+1)\eta_{2,b-1}\sha \delta_{n,b+2}\right]b+\sum\limits_{k=4}^{b}(-1)^{n+k}\left(k\eta_{2,k-1}\sha \delta_{n,k+1}\right).$$

\item\label{lem:item:j=2b=1a=2}
If $j=2$, $b=1$ and $a=2$, then $S(\sigma_{A})_j^*=(-1)^{n}(1\sha \delta_{n,3})$.

\item\label{lem:item:j=2b=1ageq3}
If $j=2$, $b=1$ and $a\geq3$, then $S(\sigma_{A})_j^*=(-1)^{n}\left(1\sha \delta_{n,3}+31\sha \delta_{n,4}\right)$.

\item\label{lem:item:j=2b=2}
If $j=2$ and $b=2$, then $S(\sigma_{A})_j^*=(-1)^{n+1}\left[1\sha (3\sha\delta_{n,5})4\right]$.

\item\label{lem:item:j=2bgeq3}
If $j=2$ and $b\geq3$, then $S(\sigma_{A})_j^*=(-1)^{n+1}\left[1\sha(3\sha \delta_{n,5})4+431\sha \delta_{n,5}\right]$.

\item\label{lem:item:3leqjleqb-1} If $3\leq j\leq b-1$, then
     \begin{align*}
S(\sigma_{A})_j^*=
(-1)^{n+j+1}&\left[\left(1\sha (j-1)\eta_{2,j-2}-\delta_{j-1,j-2}\eta_{1,j-3}\right)\sha \delta_{n,j+1}+1\sha (j+1)\eta_{2,j-1}\sha \delta_{n,j+2}\right.\\
&\hspace{75mm}\left.
+\delta_{j+2,j+1}\eta_{1,j-1}\sha \delta_{n,j+3}\right].
\end{align*}

\item\label{lem:item:3leqj=b} If $3\leq j=b$, then
      \begin{align*}
S(\sigma_{A})_{j}^*=&(-1)^{n+j+1}\left[\left(1\sha (j-1)\eta_{2,j-2}-\delta_{j-1,j-2}\eta_{1,j-3}\right)\sha \delta_{n,j+1}+1\sha (j+1)\eta_{2,j-1}\sha \delta_{n,j+2}\right].
 \end{align*}

\item\label{lem:item:3leqj=b+1a=b+1} If $3\leq j=b+1$ and $a=b+1$, then
$$S(\sigma_{A})_{j}^*=(-1)^{n+j+1}\left[\left(1\sha (j-1)\eta_{2,j-2}-\delta_{j-1,j-2}\eta_{1,j-3}\right)\sha \delta_{n,j+1}\right].$$

\item\label{lem:item:3leqj=b+1ageqb+1} If $3\leq j=b+1$ and $a\geq b+2$, then
$$S(\sigma_{A})_{j}^*=(-1)^{n+j+1}\left[\left(1\sha (j-1)\eta_{2,j-2}-\delta_{j-1,j-2}\eta_{1,j-3}\right)\sha \delta_{n,j+1}-(j+1)(j-1)\eta_{1,j-2}\sha \delta_{n,j+2}\right].$$

\item\label{lem:item:b+2leqjleqnjleqa-1} If $b+2\leq j\leq n$ and $j\leq a-1$, then
$$S(\sigma_{A})_{j}^*=(-1)^{n+j}\left[(j-1)b \eta_{1,b-1}\eta_{b+1,j-2}\sha \delta_{n,j+1}+(j+1)b\eta_{1,b-1}\eta_{b+1,j-1}\sha \delta_{n,j+2}\right].$$

\item\label{lem:item:b+2leqjleqnj=a} If $b+2\leq j\leq n$ and $j=a$, then
$S(\sigma_{A})_{j}^*=(-1)^{n+j}\left[(j-1)b \eta_{1,b-1}\eta_{b+1,j-2}\sha \delta_{n,j+1}\right]$.

\item\label{lem:item:b+2leqjleqnjgeqa+1} If $b+2\leq j\leq n$ and $j\geq a+1$, then
$S(\sigma_{A})_{j}^*=0$.
\end{enumerate}
In particular, $S(\sigma_{A})_{j}^*$ is independent of $a$ for any $j$ with $1\leq j\leq b$.
\end{theorem}

\begin{remark}

Since
\begin{align*}
1\sha (j-1)\eta_{2,j-2}-\delta_{j-1,j-2}\eta_{1,j-3}=\begin{cases}
1\sha 2-21& {\rm if}\ j=3,\\
1\sha 32-321& {\rm if}\ j=4,
\end{cases}
\end{align*}
which are not cancellation-free, so the formulas given by Theorem \ref{maintheorem}\eqref{lem:item:3leqjleqb-1}, \eqref{lem:item:3leqj=b}, \eqref{lem:item:3leqj=b+1a=b+1} and \eqref{lem:item:3leqj=b+1ageqb+1} are not cancellation-free if $j\in\{3,4\}$.
For example, for the case $j\in\{3,4\}$, the corresponding cancellation-free forms of Theorem \ref{maintheorem}\eqref{lem:item:3leqjleqb-1} can be stated as follows:
\begin{enumerate}
\item[$(h')$]\label{lem:item:3leqjleqb-12} If $j\in\{3,4\}$ with $j\leq b-1$, then
     \begin{align*}
S(\sigma_{A})_j^*=
\begin{cases}
(-1)^n\left(12\sha \delta_{n,4}+1\sha 42\sha \delta_{n,5}+5412\sha \delta_{n,6}\right)&{\rm if}\ j=3,\\
(-1)^{n+1}\left[\left(132+312\right)\sha \delta_{n,5}+1\sha523\sha\delta_{n,6}+65123\sha \delta_{n,7}\right]&{\rm if}\ j=4.
\end{cases}
\end{align*}
\end{enumerate}
\delete{Indeed, for the case $j\in\{3,4\}$, the corresponding cancellation-free forms of Theorem \ref{maintheorem}\eqref{lem:item:3leqjleqb-1}, \eqref{lem:item:3leqj=b} and \eqref{lem:item:3leqj=b+1a=b+1} can be stated, respectively, as follows:
\begin{enumerate}
\item[$(h')$]\label{lem:item:3leqjleqb-12} If $3\leq j\leq b-1$, then
     \begin{align*}
S(\sigma_{A})_j^*=
\begin{cases}
(-1)^n\left(12\sha \delta_{n,4}+1\sha 42\sha \delta_{n,5}+5412\sha \delta_{n,6}\right)&{\rm if}\ j=3,\\
(-1)^{n+1}\left[\left(132+312\right)\sha \delta_{n,5}+1\sha523\sha\delta_{n,6}+65123\sha \delta_{n,7}\right]&{\rm if}\ j=4.
\end{cases}
\end{align*}
\item[$(i')$]\label{lem:item:3leqj=b2} If $3\leq j=b$, then
      \begin{align*}
S(\sigma_{A})_{j}^*
=\begin{cases}
(-1)^n\left(12\sha \delta_{n,4}+1\sha 42\sha \delta_{n,5}\right)&{\rm if}\ j=3,\\
(-1)^{n+1}\left[\left(132+312\right)\sha \delta_{n,5}+1\sha523\sha\delta_{n,6}\right]&{\rm if}\ j=4.
\end{cases}
 \end{align*}
\item[$(j')$]\label{lem:item:3leqj=b+1a=b+12} If $3\leq j=b+1$ and $a=b+1$, then
\begin{align*}
S(\sigma_{A})_{j}^*=\begin{cases}
(-1)^n\left(12\sha \delta_{n,4}\right)&{\rm if}\ j=3,\\
(-1)^{n+1}\left[\left(132+312\right)\sha \delta_{n,5}\right]&{\rm if}\ j=4.
\end{cases}
 \end{align*}
\end{enumerate}}
Moreover, by the definition of the shuffle product, the identity given by Theorem \ref{maintheorem}\eqref{lem:item:3leqj=b+1ageqb+1} is not cancellation-free for any $j\geq3$, as $\delta_{n,j+1}(j-1)\eta_{1,j-2}$ is the intersection of the shuffle sets $1\sha (j-1)\eta_{2,j-2}\sha \delta_{n,j+1}$ and $(j+1)(j-1)\eta_{1,j-2}\sha \delta_{n,j+2}$, which have opposite signs. However, we
write these identities in the above uniforms in Theorem \ref{maintheorem}, because cancellation-free descriptions will obscure the elegance and simplicity of them.
\end{remark}

The results of Theorem \ref{maintheorem} are organized according to the relationship between $j$, $b$ and $a$. For ease of understanding, we list them in Table $1$, where $b$ and $j$ are less than $6$.
\begin{center}
\ \ Table $1$: $S(\sigma_{A})_{j}^*$ in Theorem \ref{maintheorem}, where $j,b\in[5]$.\\
\vspace{2mm}
\begin{tabular}{|c|c|c|c|c|c|}
\hline
$S(\sigma_A)_j^*$&$j=1$&$j=2$&$j=3$&$j=4$&$j=5$\\
\hline
$b=1$&\eqref{lem:item:j=1b=1}
&\makecell{\eqref{lem:item:j=2b=1a=2} if $a=2$,\\ \eqref{lem:item:j=2b=1ageq3} if $a\geq3$.}
&\makecell{\eqref{lem:item:b+2leqjleqnjgeqa+1} if $a=2$,\\ \eqref{lem:item:b+2leqjleqnj=a} if $a=3$,\\ \eqref{lem:item:b+2leqjleqnjleqa-1} if $a\geq4$.}
&\makecell{\eqref{lem:item:b+2leqjleqnjgeqa+1} if $2\leq a\leq 3$,\\ \eqref{lem:item:b+2leqjleqnj=a} if $a=4$,\\ \eqref{lem:item:b+2leqjleqnjleqa-1} if $a\geq5$.}
&\makecell{\eqref{lem:item:b+2leqjleqnjgeqa+1} if $2\leq a\leq 4$,\\ \eqref{lem:item:b+2leqjleqnj=a} if $a=5$,\\ \eqref{lem:item:b+2leqjleqnjleqa-1} if $a\geq6$.}\\
\hline
$b=2$&\eqref{lem:item:j=1b=2}
&\eqref{lem:item:j=2b=2}
&\makecell{\eqref{lem:item:3leqj=b+1a=b+1} if $a=3$,\\ \eqref{lem:item:3leqj=b+1ageqb+1} if $a\geq4$.}
&\makecell{\eqref{lem:item:b+2leqjleqnjgeqa+1} if $a=3$,\\ \eqref{lem:item:b+2leqjleqnj=a} if $a=4$,\\ \eqref{lem:item:b+2leqjleqnjleqa-1} if $a\geq5$.}
&\makecell{\eqref{lem:item:b+2leqjleqnjgeqa+1} if $3\leq a\leq 4$,\\ \eqref{lem:item:b+2leqjleqnj=a} if $a=5$,\\ \eqref{lem:item:b+2leqjleqnjleqa-1} if $a\geq6$.}\\
\hline
$b=3$&\eqref{lem:item:j=1bgeq3}
&\eqref{lem:item:j=2bgeq3}
&\eqref{lem:item:3leqj=b}&\makecell{\eqref{lem:item:3leqj=b+1a=b+1} if $a=4$,\\ \eqref{lem:item:3leqj=b+1ageqb+1} if $a\geq5$.}
&\makecell{\eqref{lem:item:b+2leqjleqnjgeqa+1} if $a=4$,\\ \eqref{lem:item:b+2leqjleqnj=a} if $a=5$,\\ \eqref{lem:item:b+2leqjleqnjleqa-1} if $a\geq6$.} \\
\hline
$b=4$&\eqref{lem:item:j=1bgeq3}
&\eqref{lem:item:j=2bgeq3}
&\eqref{lem:item:3leqjleqb-1}
&\eqref{lem:item:3leqj=b}
&\makecell{\eqref{lem:item:3leqj=b+1a=b+1} if $a=5$,\\ \eqref{lem:item:3leqj=b+1ageqb+1} if $a\geq6$.}\\
\hline
$b=5$&\eqref{lem:item:j=1bgeq3}
&\eqref{lem:item:j=2bgeq3}
&\eqref{lem:item:3leqjleqb-1}
&\eqref{lem:item:3leqjleqb-1}
&\eqref{lem:item:3leqj=b}\\
\hline
\end{tabular}
\end{center}

\vspace{5mm}
\begin{exmp}\label{exmp234n}
For $n\in\{2,3,4\}$ and $A\subseteq[n]$ with $|A|=2$, applying Eq.\,\eqref{eq:antipodeformula1}, a simple recursive computation yields the formulas for $S(\sigma_{A,n})$, where $S(\sigma_{A,n})_j^*$ are listed in Table $2$.
\begin{center}
\ \ Table $2$: $S(\sigma_{A,n})_{j}^*$ for $n\in\{2,3,4\}$.\\
\vspace{2mm}
\begin{tabular}{|c|c|c|c|c|}
\hline
$S(\sigma_A)_j^*$&$j=1$&$j=2$&$j=3$&$j=4$\\
\hline
$\sigma_{\{1,2\},2}$&$0$&$1$&$0$&$0$\\
\hline
$\sigma_{\{1,2\},3}$&$23$&$-(1\sha3)$&$0$&$0$\\
\hline
$\sigma_{\{1,3\},3}$&$23$&$-\left(1\sha3+31\right)$&$21$&$0$ \\
\hline
$\sigma_{\{2,3\},3}$&$0$&$0$&$-12$&$0$\\
\hline
$\sigma_{\{1,2\},4}$&$-(2\sha4)3$&$1\sha43$&$0$&$0$\\
\hline
$\sigma_{\{1,3\},4}$&$-(2\sha4)3$&$1\sha43+31\sha4$&$-(21\sha 4)$&$0$\\
\hline
$\sigma_{\{1,4\},4}$&$-(2\sha4)3$&$1\sha43+31\sha4$&$-(21\sha4+412)$&$312$\\
\hline
$\sigma_{\{2,3\},4}$&$324$&$-(1\sha34)$&$(1\sha2-21)\sha4$&$0$\\
\hline
$\sigma_{\{2,4\},4}$&$324$&$-(1\sha 34)$&$(1\sha2-21)\sha4-421$&$321$\\
\hline
$\sigma_{\{3,4\},4}$&$423$&$-(1\sha 34+431)$&$(1\sha2-21)\sha4+1\sha42$&$-(1\sha32-321)$\\
\hline
\end{tabular}
\end{center}
\vspace{3mm}
As an illuminating illustration, we give the details of the computations of $S(4312)$, the other formulas for $S(\sigma_{A,n})$ can be similarly obtained.
By \cite[Propositions 9.3 and 9.5]{BS17}, or by Lemma \ref{lem:lema1} below, we have $S(1)=-1$, $S(21)=12$ and $S(321)=-123$.
Then
\delete{\begin{align*}
S(312)&=-312-S(\st(3))\overline{\sha} \st(12)-S(\st(31))\overline{\sha}\st(2)&{\rm by\ Eq.\,\eqref{eq:antipodeformula1}}\\
&=-312-S(1)\sha 23-S(21)\sha 3\\
&=-312+1\sha 23-12\sha3\\
&=-312+231+(1\sha 2)3-(1\sha3)2-123&{\rm by\ Eq.\,\eqref{eq:uashavb}}
\end{align*}
and}
\begin{align*}
S(4312)&=-4312-S(\st(4))\overline{\sha} \st(312)-S(\st(43))\overline{\sha}\st(12)-S(\st(431))\overline{\sha}\st(2)&{\rm by\ Eq.\,\eqref{eq:antipodeformula1}}\\
&=-4312-S(1)\sha 423-S(21)\sha34-S(321)\sha4\\
&=-4312+1\sha 423-12\sha34+123\sha4\\
&=-4312+4231+(1\sha 42)3-(1\sha34)2-(12\sha3)4+1234+(12\sha4)3&{\rm by\ Eq.\,\eqref{eq:uashavb}}.
\end{align*}
By grouping together terms ending in the same characters, we obtain that
\begin{align*}
%S(312)&=231-(1\sha3+31)2+(1\sha 2-12)3=231-(1\sha3+31)2+213,\\
S(4312)&=4231-(1\sha 34+431)2+\left(1\sha42+12\sha4\right)3-(12\sha3-123)4\\
&=4231-(1\sha 34+431)2+\left[(1\sha2-21)\sha4+1\sha42\right]3-(1\sha32-321)4.
\end{align*}
Hence,
\begin{align*}
S(4312)_1^*&=423,& S(4312)_2^*&=-(1\sha 34+431),\\
S(4312)_3^*&=(1\sha2-21)\sha4+1\sha42,&S(4312)_4^*&=-(1\sha32-321).
\end{align*}
These computations can be checked by the SageMath computer algebra system \cite{SG23}, where packages on the Malvenuto-Reutenauer Hopf algebra \cite{CG17}
and on words shuffles \cite{Pr17} are helpful. For example, if we write $4312=F_{4312}$, then $S(4312)$ can be verified by the following code:
\begin{lstlisting}
sage: FQSym = algebras.FQSym(QQ)
sage: F = FQSym.F()
sage: F[[4,3,1,2]].antipode()
\end{lstlisting}
\delete{Completely analogous to the previous computations, using Eqs.\,\eqref{eq:uashavb} and \eqref{eq:antipodeformula1}, we obtain that
\begin{align*}
S(213)&%=-213-S(1)\sha 23-S(21)\sha 3=-213+1\sha 23-12\sha 3
=231-(1\sha3)2,\\
S(312)&=231-(1\sha3+31)2+213,\\
S(2134)%&=-2134-S(1)\sha 234-S(21)\sha 34-S(213)\sha4\\
%&=-2134+1\sha 234-12\sha 34-231\sha4+(1\sha3)2\sha4\\
&=-(2\sha 4)31+(1\sha43)2,\\
S(3124)%&=-3124-S(1)\sha 234-S(21)\sha 34-S(312)\sha4\\
%&=-3124+1\sha 234-12\sha 34-[231-(1\sha3+31)2+213]\sha4\\
&=-(2\sha 4)31+(1\sha43+31\sha4)2-(21\sha4)3,\\
S(4123)%&=-4123-S(1)\sha 234-S(21)\sha 34-S(312)\sha4\\
%&=-4123+1\sha 234-12\sha 34-[231-(1\sha3+31)2+213]\sha4\\
&=-(2\sha 4)31+(1\sha43+31\sha4)2-(21\sha4+412)3+3124,\\
S(3214)%&=-3214-S(1)\sha 324-S(21)\sha 34-S(321)\sha4\\
%&=-3214+1\sha 324-12\sha 34+123\sha4\\
&=3241-(1\sha 34)2+[(1\sha2-21)\sha4]3,\\
S(4213)%&=-4213-S(1)\sha 324-S(21)\sha 34-S(321)\sha4\\
%&=-4213+1\sha 324-12\sha 34+123\sha4\\
&=3241-(1\sha 34)2+[(1\sha2-21)\sha4-421]3+3214.
\end{align*}}
\end{exmp}

\section{Proof of the main theorem}\label{Proofofthemainythm}
In this section, we prove Theorem \ref{maintheorem}, from which \cite[Conjectures 9.9 and 9.10]{BS17} follows easily.
We begin with some useful lemmas, where the first one was proved in \cite{BS17}.

\begin{lemma}\cite[Propositions 9.3 and 9.5]{BS17}\label{lem:lema1}
For any $n\in\mathbb{P}$, we have
\begin{equation*}
S(\eta_{1,n}) = (-1)^{n}\delta_{n,1}\quad \text{and}\quad  S(\delta_{n,1}) = (-1)^{n}\eta_{1,n}.
\end{equation*}
\end{lemma}

\begin{lemma}\label{lem:-1ideltaij+1shaeta}
 For any $k\in\mathbb{N}$ and $n\in\mathbb{P}$ with $k\leq n$, we have
\begin{align*}
\sum_{i=k}^{n-1}(-1)^{i} \delta_{i,k+1}\sha  \eta_{i+1,n}=(-1)^{n+1} \delta_{n,k+1}.
\end{align*}
\end{lemma}
\begin{proof} For each $i$ with $k+1\leq i\leq n-1$, it follows from Eq.\eqref{eq:aushabv} that
\begin{align*}
\delta_{i,k+1}\sha  \eta_{i+1,n}=i\left(\delta_{i-1,k+1}\sha\eta_{i+1,n}\right)+(i+1)\left(\delta_{i,k+1}\sha\eta_{i+2,n}\right).
\end{align*}
Thus,
\begin{align*}
\sum_{i=k}^{n-1}(-1)^{i} \delta_{i,k+1}\sha  \eta_{i+1,n}
&=(-1)^k\eta_{k+1,n}+\sum_{i=k+1}^{n-1}(-1)^{i}i\left(\delta_{i-1,k+1}\sha\eta_{i+1,n}\right)+\sum_{i=k+1}^{n-1}(-1)^{i}(i+1)\left(\delta_{i,k+1}\sha\eta_{i+2,n}\right).
\end{align*}
Since $i\left(\delta_{i-1,k+1}\sha\eta_{i+1,n}\right)=\eta_{k+1,n}$ if $i=k+1$, and
$(i+1)\left(\delta_{i,k+1}\sha\eta_{i+2,n}\right)=\delta_{n,k+1}$ if $i=n-1$,
the proof follows by substituting $i$ by $i+1$ in the second summand in the above identity.
\end{proof}

\begin{lemma}\label{lem:eta2n+-1i}
For any positive integers $j,k$ and $n$ with $j\leq k\leq n-1$, we have
\begin{equation}\label{etak+1n1}
(-1)^k\delta_{k,j}\eta_{k+1,n} + \sum_{i=k+1}^{n-1}(-1)^{i}\left(\delta_{k,j} \sha \delta_{i,k+2}\right)(k+1) \sha \eta_{i+1,n}=(-1)^{n+1}\left(\delta_{k,j} \sha \delta_{n,k+2}\right)(k+1).
\end{equation}
\end{lemma}
\begin{proof}
The proof is by induction on $n$. If $n=k+1$, then the second term in the left-hand side of Eq.\,\eqref{etak+1n1} is $0$, so that Eq.\,\eqref{etak+1n1} is true.
Now assume that Eq.\,\eqref{etak+1n1} holds for $n-1$. For $k+1\leq i\leq n-1$, using Eq.\,\eqref{eq:uashavb}, we get
\begin{align*}
\left(\delta_{k,j} \sha \delta_{i,k+2}\right)(k+1) \sha \eta_{i+1,n}=\left(\delta_{k,j} \sha \delta_{i,k+2} \sha \eta_{i+1,n}\right)(k+1)
+\left[\left(\delta_{k,j} \sha \delta_{i,k+2}\right)(k+1) \sha \eta_{i+1,n-1}\right]n.
\end{align*}
Thus, by grouping the shuffle sets according to their last characters, the left-hand side of Eq.\,\eqref{etak+1n1} is equal to
\begin{align*}
&{}\sum_{i=k+1}^{n-1}(-1)^{i}\left(\delta_{k,j} \sha \delta_{i,k+2} \sha \eta_{i+1,n}\right)(k+1)\\
&{}\hspace{30mm}+\left((-1)^k\delta_{k,j}\eta_{k+1,n-1} +\sum_{i=k+1}^{n-1}(-1)^{i}\left(\delta_{k,j} \sha \delta_{i,k+2}\right)(k+1) \sha \eta_{i+1,n-1}\right)n\\
=&\,\sum_{i=k+1}^{n-1}(-1)^{i}\left(\delta_{k,j} \sha \delta_{i,k+2} \sha \eta_{i+1,n}\right)(k+1)\\
=&\,\left[\delta_{k,j} \sha \left(\sum_{i=k+1}^{n-1}(-1)^{i}\delta_{i,k+2} \sha \eta_{i+1,n}\right)\right](k+1)\\
=&\,(-1)^{n+1}(\delta_{k,j} \sha \delta_{n,k+2})(k+1),
\end{align*}
where the first equality follows from the induction hypothesis, while the third equality follows from Lemma \ref{lem:-1ideltaij+1shaeta}.
\end{proof}

The proof of the following identity of shuffle sets is analogous to that of Lemma \ref{lem:eta2n+-1i}, so we will just sketch the basic facts and omit most
details.
\begin{lemma}\label{lem:etadeltaetan}
For any positive integers $k$ and $n$ with $2\leq k\leq n-1$, we have
\begin{align*}
    (-1)^{k}(k+1)\eta_{2,k}\eta_{k+2,n}+\sum_{i=k+1}^{n-1}(-1)^{i}\left[(k+1)\eta_{2,k-1}\sha \delta_{i,k+2}\right]k\sha\eta_{i+1,n}
    =(-1)^{n+1}\left[(k+1)\eta_{2,k-1}\sha \delta_{n,k+2}\right]k.
\end{align*}
\end{lemma}
\begin{proof}
We act again by induction on $n$, with the base case $n=k+1$ being trivial. Now assume that the desired identity holds for $n-1$.
According to Eq.\,\eqref{eq:uashavb}, for any $i$ with $k+1\leq i\leq n-1$, we have
\begin{align*}
\left[(k+1)\eta_{2,k-1}\sha \delta_{i,k+2}\right]k\sha\eta_{i+1,n}
=&\left[(k+1)\eta_{2,k-1}\sha \delta_{i,k+2}\sha\eta_{i+1,n}\right]k\\
&{}\hspace{15mm}+\left[\left((k+1)\eta_{2,k-1}\sha \delta_{i,k+2}\right)k\sha\eta_{i+1,n-1}\right]n.
\end{align*}
The remainder of the proof is analogous to that of Lemma \ref{lem:eta2n+-1i}; the details are omitted.
\end{proof}
\begin{lemma}\label{lem:etai-n}
Let $i$ and $n$ be positive integers with $2\leq i\leq n$, and let $\sigma=\sigma_1\cdots \sigma_{n}\in\mathfrak{S}_n$ such that $\sigma_i\cdots \sigma_{n}=\eta_{i,n}$. Then for any positive integer $j$ with $i\leq j\leq n$, we have $S(\sigma)_j=0$.
\end{lemma}
\begin{proof}
We proceed by induction on $n$, the base case $n=2$ is trivial by Lemma \ref{lem:lema1}.
Assume that the statement is true for all positive integers less than $n$, where $n\geq3$.
If $i\leq j\leq n-1$,
then it follows from Eq.\,\eqref{eq:antipodeformula1} that all shuffles ending in the character $j$ must be contained in
\begin{align*}
-\sum_{k=j}^{n-1}S(\st(\sigma_1\cdots \sigma_{i-1}\eta_{i,k}))\overline{\sha}\st(\eta_{k+1,n})
=-\sum_{k=j}^{n-1}S(\sigma_1\cdots \sigma_{i-1}\eta_{i,k})\sha\eta_{k+1,n}.
\end{align*}
So $S(\sigma)_j$ is only contained in the shuffle sets $S(\sigma_1\cdots \sigma_{i-1}\eta_{i,k})_j\sha\eta_{k+1,n}$, where $j\leq k\leq n-1$.
By the induction hypothesis,  $S(\sigma_1\cdots \sigma_{i-1}\eta_{i,k})_j=0$ for all $k$ with $j\leq k\leq n-1$. Thus, $S(\sigma)_j=0$ for any $j$ with $i\leq j\leq n-1$.
Consider now the case $j=n$.
Note that $\sigma_n=n$, it follows from Eq.\,\eqref{eq:antipodeformula1} that
\begin{align*}
S(\sigma)_n^*&=-\sum_{k=0}^{n-1}S(\sigma_1\cdots \sigma_{k})\overline{\sha}\st(\sigma_{k+1}\cdots \sigma_{n-1})\\
&=-S(\sigma_1\cdots \sigma_{n-1})-\sum_{k=0}^{n-2}S(\sigma_1\cdots \sigma_{k})\overline{\sha}\st(\sigma_{k+1}\cdots \sigma_{n-1}).
\end{align*}
Again by Eq.\,\eqref{eq:antipodeformula1}, we have
\begin{align*}
S(\sigma_1\cdots \sigma_{n-1})=-\sum_{k=0}^{n-2}S(\sigma_1\cdots \sigma_{k})\overline{\sha}\st(\sigma_{k+1}\cdots \sigma_{n-1}),
\end{align*}
so that $S(\sigma)_n^*=0$. Hence, $S(\sigma)_n=0$, completing the proof.
\end{proof}

We are now in a position to prove Theorem \ref{maintheorem}. The proof is by induction on $n$, where the base cases $n\in\{2,3,4\}$
have been checked by Example \ref{exmp234n}. We first apply the antipode formula given by Eq.\,\eqref{eq:antipodeformula1} to $\sigma_A$
to obtain the general form of $S(\sigma_A)$ (i.e., Eq.\,\eqref{eq:s(ba--n)}), and then proceed by a careful analysis of the relationship between $j$ and $a,b$.

To begin with, we specialize Theorem \ref{maintheorem} to the case $a=n$ and $b=n-1$, which is a special case need to be analyzed separately in the proof of Theorem \ref{maintheorem}.

Observe that for any words $u,v$ on an alphabet $A$, and any element $j$ in $A$, we have
\begin{align}\label{eq:ushavjeq}
(u\sha v)_j^*=u_j^*\sha v+u\sha v_j^*.
\end{align}
In the proofs of Lemma \ref{lem:etaibn-1} and Theorem \ref{maintheorem},
we make frequent use of the above observation.% but not mentioned in the corresponding applications.

\begin{lemma}\label{lem:etaibn-1}
Let $A=\{n-1,n\}$, where $n\in\mathbb{P}$ with $n\geq4$.
\begin{enumerate}
  \item\label{lem:item:sba1a-1'}
$S(\sigma_{A})_1^*=n\eta_{2,n-1}+\sum\limits_{k=4}^{n-1}(-1)^{n+k}\left(k\eta_{2,k-1}\sha \delta_{n,k+1}\right)$.
  \item\label{lem:item:sba1a-1j=2'}
$S(\sigma_{A})_2^*=(-1)^{n+1}\left[1\sha(3\sha \delta_{n,5})4+431\sha \delta_{n,5}\right].$

   \item\label{lem:item:sbaa'} For $3\leq j\leq n-2$, we have
     \begin{align*}
S(\sigma_{A})_j^*=
(-1)^{n+j+1}&\left[\left(1\sha (j-1)\eta_{2,j-2}-\delta_{j-1,j-2}\eta_{1,j-3}\right)\sha \delta_{n,j+1}+1\sha (j+1)\eta_{2,j-1}\sha \delta_{n,j+2}\right.\\
 &\hspace{73mm}+\left.\delta_{j+2,j+1}\eta_{1,j-1}\sha \delta_{n,j+3}\right].
\end{align*}

   \item\label{lem:item:sbaaj=b'}
$S(\sigma_{A})_{n-1}^*=1\sha n\eta_{2,n-2}+\left[1\sha (n-2)\eta_{2,n-3}-\delta_{n-2,n-3}\eta_{1,n-4}\right]\sha n.$

   \item\label{lem:item:sbaaj=b+1geq3'} $S(\sigma_{A})_{n}^*=\delta_{n-1,n-2}\eta_{1,n-3}-1\sha (n-1)\eta_{2,n-2}$.
\end{enumerate}
\end{lemma}
\begin{proof}
The proof is by induction on $n$.
In view of Example \ref{exmp234n}, the assertion is true for $n=4$. We assume that parts \eqref{lem:item:sba1a-1'}-\eqref{lem:item:sbaaj=b+1geq3'} of this lemma are true for any positive integer less than $n$, where $n\geq5$.
By Lemma \ref{lem:lema1}, $S(1) =-1$, $S(21) = 12$, $S(321) = -123$. Since $\sigma_A=\delta_{n,n-1}\eta_{1,n-2}$, it follows from Eq.\,\eqref{eq:antipodeformula1} that
\begin{align}\label{eq:s(ba--ngeq5n-1n)}
 S(\sigma_A)_j^*=U_j^*+V_j^*+W_j^*+X_j^*+Y_j^*,
\end{align}
where
\begin{align*}
 U=-\delta_{n,n-1}\eta_{1,n-2},\  V=1\sha n\eta_{2,n-1},\  W=-12\sha \eta_{3,n},\  X=123\sha \eta_{4,n},\  Y=-\sum_{i=4}^{n-1} S(\sigma_{\{i-1,i\},i})\sha\eta_{i+1,n}.
\end{align*}
We now determine $S(\sigma_A)_j^*$  according to the values of $j$.

\eqref{lem:item:sba1a-1'} By Eqs.\,\eqref{eq:ushavjeq} and \eqref{eq:s(ba--ngeq5n-1n)}, together with the induction hypothesis on part \eqref{lem:item:sba1a-1'}, we have
\begin{align*}
 S(\sigma_A)_1^*&=V_1^*+Y_1^*\\
 &= n\eta_{2,n-1} -\sum_{i=4}^{n-1} S(\sigma_{\{i-1,i\},i})_1^*\sha\eta_{i+1,n}\\
 &=n\eta_{2,n-1} -\sum_{i=4}^{n-1}i\eta_{2,i-1}\sha\eta_{i+1,n}+\sum_{i=4}^{n-1}\sum_{k=4}^{i-1}(-1)^{i+k+1}k\eta_{2,k-1}\sha\delta_{i,k+1}\sha\eta_{i+1,n}.
\end{align*}
By Lemma \ref{lem:-1ideltaij+1shaeta}, the third summand in the above equation is
\begin{align*}
\sum_{k=4}^{n-2}(-1)^{k+1}k\eta_{2,k-1}\sha\left(\sum_{i=k+1}^{n-1}(-1)^{i}\delta_{i,k+1}\sha\eta_{i+1,n}\right)
=&\,\sum_{k=4}^{n-2}(-1)^{k+1}k\eta_{2,k-1}\sha\left((-1)^{k+1}\eta_{k+1,n}+(-1)^{n+1}\delta_{n,k+1}\right)\\
=&\,\sum_{k=4}^{n-2}k\eta_{2,k-1}\sha\eta_{k+1,n}+\sum_{k=4}^{n-2}(-1)^{n+k}k\eta_{2,k-1}\sha\delta_{n,k+1}.
\end{align*}
Hence,
\begin{align*}
 S(\sigma_A)_1^* &=n\eta_{2,n-1} -(n-1)\eta_{2,n-2}\sha n+\sum_{k=4}^{n-2}(-1)^{n+k}k\eta_{2,k-1}\sha\delta_{n,k+1}\\
 &=n\eta_{2,n-1}+\sum_{k=4}^{n-1}(-1)^{n+k}k\eta_{2,k-1}\sha\delta_{n,k+1},
\end{align*}
as required.

\eqref{lem:item:sba1a-1j=2'}
Applying Eqs.\,\eqref{eq:ushavjeq} and \eqref{eq:s(ba--ngeq5n-1n)}, together with the induction hypothesis on part \eqref{lem:item:sba1a-1j=2'}, yields that
\begin{align*}
 S(\sigma_A)_2^*&=W_2^*+Y_2^*\\
 &=-1\sha \eta_{3,n}-\sum_{i=4}^{n-1} S(\sigma_{\{i-1,i\},i})_2^*\sha\eta_{i+1,n}\\
 &=-1\sha \eta_{3,n}+\sum_{i=4}^{n-1} (-1)^i\left[1\sha(3\sha \delta_{i,5})4+431\sha \delta_{i,5}\right]\sha\eta_{i+1,n}\\
 &=1\sha \left[-\eta_{3,n}+\sum_{i=4}^{n-1} (-1)^i(3\sha \delta_{i,5})4\sha\eta_{i+1,n}\right]+431\sha\left(\sum_{i=4}^{n-1} (-1)^i \delta_{i,5}\sha\eta_{i+1,n}\right).
\end{align*}
Since $n\geq4$, by Lemma \ref{lem:eta2n+-1i} with $k=j=3$ and Lemma \ref{lem:-1ideltaij+1shaeta} with $k=4$, we have
\begin{align*}
 S(\sigma_A)_2^*=1\sha\left((-1)^{n+1}(3\sha \delta_{n,5})4\right)+431\sha\left((-1)^{n+1}\delta_{n,5}\right)=(-1)^{n+1}\left[1\sha(3\sha \delta_{n,5})4+431\sha \delta_{n,5}\right].
\end{align*}

\eqref{lem:item:sbaa'}
We first consider the case where $3\leq j=n-2$. Then, by Eqs.\,\eqref{eq:ushavjeq} and \eqref{eq:s(ba--ngeq5n-1n)}, we get
$S(\sigma_A)_j^*=U_j^*+X_j^*+Y_j^*$, where $X_j^*$ and $Y_j^*$ can be merged, so that
\begin{align*}%\label{eq:s(ba--n)3leqjleqb-1}
 S(\sigma_A)_j^*=-\delta_{n,n-1}\eta_{1,n-3}-\sum_{i=n-2}^{n-1} S(\sigma_{\{i-1,i\},i})_{n-2}^*\sha\eta_{i+1,n}.
\end{align*}
According to the induction hypothesis, that is, by part \eqref{lem:item:sbaaj=b+1geq3'} if $i=n-2$, and by part \eqref{lem:item:sbaaj=b'} if $i=n-1$, we obtain that
\begin{align*}%\label{eq:s(ba--n)3leqjleqb-1}
S(\sigma_{\{i-1,i\},i})_{n-2}^*=\begin{cases}
\delta_{n-3,n-4}\eta_{1,n-5}-1\sha (n-3)\eta_{2,n-4}&{\rm if}\ i=n-2,\\
1\sha (n-1)\eta_{2,n-3}+\left[1\sha (n-3)\eta_{2,n-4}-\delta_{n-3,n-4}\eta_{1,n-5}\right]\sha(n-1)&{\rm if}\ i=n-1.
\end{cases}
\end{align*}
So
\begin{align*}%\label{eq:s(ba--n)3leqjleqb-1}
 S(\sigma_A)_j^*
 =&-\delta_{n,n-1}\eta_{1,n-3}-\left(\delta_{n-3,n-4}\eta_{1,n-5}-1\sha (n-3)\eta_{2,n-4}\right)\sha\eta_{n-1,n}\\
 &-\left[1\sha (n-1)\eta_{2,n-3}\sha n+\left(1\sha (n-3)\eta_{2,n-4}-\delta_{n-3,n-4}\eta_{1,n-5}\right)\sha(n-1)\sha n\right]\\
 =&-\delta_{n,n-1}\eta_{1,n-3}-\left[1\sha (n-3)\eta_{2,n-4}-\delta_{n-3,n-4}\eta_{1,n-5}\right]\sha\delta_{n,n-1}-1\sha (n-1)\eta_{2,n-3}\sha n,
\end{align*}
agreeing with the identity given by part \eqref{lem:item:sbaa'}.

 We next consider the case where $3\leq j\leq n-3$. It follows from Eqs.\,\eqref{eq:ushavjeq} and \eqref{eq:s(ba--ngeq5n-1n)} that
 \begin{align}\label{eq:s(ba--n)3leqjleqb-111}
 S(\sigma_A)_j^*=X_j^*+Y_j^*=-\sum_{i=j}^{n-1} S(\sigma_{\{i-1,i\},i})_j^*\sha\eta_{i+1,n}.
\end{align}
Applying part \eqref{lem:item:sbaaj=b+1geq3'} if $i=j$, part \eqref{lem:item:sbaaj=b'} if $i=j+1$, and part \eqref{lem:item:sbaa'} if $j+2\leq i\leq n-1$, to $S(\sigma_{\{i-1,i\},i})_j^*$,
we see from the induction hypothesis that
\begin{align*}
S(\sigma_A)_j^*=&\,-\left[\delta_{j-1,j-2}\eta_{1,j-3}-1\sha (j-1)\eta_{2,j-2}\right]\sha\eta_{j+1,n}\\
&-\left[1\sha (j+1)\eta_{2,j-1}+\left(1\sha (j-1)\eta_{2,j-2}-\delta_{j-1,j-2}\eta_{1,j-3}\right)\sha (j+1)\right]\sha\eta_{j+2,n}\\
&\,+\sum_{i=j+2}^{n-1}(-1)^{i+j}\left[\left(1\sha (j-1)\eta_{2,j-2}-\delta_{j-1,j-2}\eta_{1,j-3}\right)\sha \delta_{i,j+1}+1\sha (j+1)\eta_{2,j-1}\sha \delta_{i,j+2}\right.\\
 &\hspace{90mm}+\left.\delta_{j+2,j+1}\eta_{1,j-1}\sha \delta_{i,j+3}\right]\sha\eta_{i+1,n}\\
 =&\,(-1)^{j}\left[\left(1\sha (j-1)\eta_{2,j-2}-\delta_{j-1,j-2}\eta_{1,j-3}\right)\sha\left( \sum_{i=j}^{n-1}(-1)^{i}\delta_{i,j+1}\sha\eta_{i+1,n}\right)\right.\\
&\hspace{-10mm}+\left.\left(1\sha (j+1)\eta_{2,j-1}\right)\sha \left(\sum_{i=j+1}^{n-1}(-1)^{i}\delta_{i,j+2}\sha\eta_{i+1,n}\right)
 +\left(\delta_{j+2,j+1}\eta_{1,j-1}\right)\sha\left(\sum_{i=j+2}^{n-1}(-1)^{i}\delta_{i,j+3}\sha\eta_{i+1,n}\right)\right].
\end{align*}
The proof now follows from Lemma \ref{lem:-1ideltaij+1shaeta} with $k\in\{j,j+1,j+2\}$.

\eqref{lem:item:sbaaj=b'} Using part \eqref{lem:item:sbaaj=b+1geq3'}, where we substitute $n$ by $n-1$, it follows from the induction hypothesis on part \eqref{lem:item:sbaaj=b'} and Eqs.\,\eqref{eq:ushavjeq} and \eqref{eq:s(ba--ngeq5n-1n)} that
\begin{align*}%\label{eq:s(ba--ngeq5n-1n)}
 S(\sigma_A)_{n-1}^*&=V_{n-1}^*+Y_{n-1}^*\\
&=1\sha n\eta_{2,n-2}-S(\sigma_{\{n-2,n-1\},n-1})_{n-1}^*\sha n\\
&=1\sha n\eta_{2,n-2}+\left(1\sha(n-2)\eta_{2,n-3}-\delta_{n-2,n-3}\eta_{1,n-4}\right)\sha n,
\end{align*}
as required.

\eqref{lem:item:sbaaj=b+1geq3'} From Eqs.\,\eqref{eq:ushavjeq} and \eqref{eq:s(ba--ngeq5n-1n)} we see that
\begin{align*}%\label{eq:s(ba--ngeq5n-1n)}
 S(\sigma_A)_n^*=W_{n}^*+X_{n}^*+Y_{n}^*%=-\sum_{i=2}^{n-1} S(\sigma_{\{i-1,i\},i})\sha\eta_{i+1,n-1}
 =-\sum_{i=2}^{n-2} S(\sigma_{\{i-1,i\},i})\sha\eta_{i+1,n-1}-S(\sigma_{\{n-2,n-1\},n-1}).
\end{align*}
Substituting $n$ by $n-1$ in Eq.\,\eqref{eq:s(ba--ngeq5n-1n)} yields that
\begin{align*}%\label{eq:s(ba--ngeq5n-1n)}
S(\sigma_{\{n-2,n-1\},n-1})=-\delta_{n-1,n-2}\eta_{1,n-3}+1\sha (n-1)\eta_{2,n-2}-\sum_{i=2}^{n-2}S(\sigma_{\{i-1,i\},i})\sha\eta_{i+1,n-1},
\end{align*}
and hence
\begin{align*}%\label{eq:s(ba--ngeq5n-1n)}
 S(\sigma_A)_n^*&=\delta_{n-1,n-2}\eta_{1,n-3}-1\sha (n-1)\eta_{2,n-2}.
\end{align*}
So the proof follows by induction.
\end{proof}

\begin{proof}[Proof of Theorem \ref{maintheorem}]
We proceed by induction on $n$. For $n\in\{2,3,4\}$, the assertion is true by Example \ref{exmp234n}.
Now assume that Theorem \ref{maintheorem} is true for any positive integer less than $n$, where $n\geq5$.
By Lemma \ref{lem:etaibn-1}, we need only assume that $b\leq n-2$.
It follows from Eq.~\eqref{eq:antipodeformula1} and Lemma \ref{lem:lema1} that $S(\sigma_A)=U+V+W+X+Y+Z$, and hence for all $j\in[n]$, we have
\begin{align}\label{eq:s(ba--n)}
 S(\sigma_A)_j^*=U_j^*+V_j^*+W_j^*+X_j^*+Y_j^*+Z_j^*,
\end{align}
where
\begin{align*}
	U&=-ab\eta_{1,b-1}\eta_{b+1,a-1}\eta_{a+1,n},
	&V&=1\sha (b+1)\eta_{2,b}\eta_{b+2,n},
    &W&=-12\sha \eta_{3,n},\\
	X&=-\sum_{i=3}^{b+1} S(\sigma_{\{i-1,i\},i})\sha\eta_{i+1,n},
	&Y&=-\sum_{i=b+2}^{{\rm min}\{a,n-1\}}S(\sigma_{\{b,i\},i})\sha\eta_{i+1,n},
	&Z&=-\sum_{i=a+1}^{n-1}S(\sigma_{\{b,a\},i})\sha\eta_{i+1,n}.
\end{align*}
Now we show that for each $j$ with $1 \leq j \leq n$, the term $S(\sigma_A)_j^*$ coincides with the desired form.

We need two simple but useful observations. First, note that for any positive integers $i$ and $j$ with $1\leq j\leq b$ and $ a\leq i\leq n-1$, the induction hypothesis guarantees that $S(\sigma_{\{b,a\},i})_{j}^*$ is independent of $a$. So we can substitute $a$ by $i$ in $S(\sigma_{\{b,a\},i})_{j}^*$  without changing its value, that is, we have $S(\sigma_{\{b,a\},i})_{j}^*=S(\sigma_{\{b,i\},i})_{j}^*$, where $1\leq j\leq b$ and $a\leq i\leq n-1$. Hence, $Y_j^*$ and $Z_j^*$ can be merged into one term, that is,
\begin{align}\label{eq:nn-1--1jb-1}
Y_j^*+Z_j^*=-\sum_{i=b+2}^{n-1}S(\sigma_{\{b,i\},i})_j^*\sha\eta_{i+1,n}, \ {\rm where}\ 1\leq j\leq b\ \text{and}\ a\leq i\leq n-1.
\end{align}
Second, if $b\geq2$ and $i=3$, then $X$ is nonzero, and we have
\begin{align}\label{eq:ageq2i=3S}
    S(\sigma_{\{i-1,i\},i})\sha\eta_{i+1,n}=S(\sigma_{\{2,3\},3})\sha\eta_{4,n}=S(321)\sha\eta_{4,n}=-123\sha\eta_{4,n}.
\end{align}

We first prove parts \eqref{lem:item:j=1b=1}, \eqref{lem:item:j=1b=2} and \eqref{lem:item:j=1bgeq3}. Assume that $j=1$. Then $U_j^*=W_j^*=0$.
According to Eq.\,\eqref{eq:ageq2i=3S}, we have $(S(\sigma_{\{i-1,i\},i})\sha\eta_{i+1,n})_1^*=0$ if $i=3$.
It follows from Eqs.\,\eqref{eq:ushavjeq}, \eqref{eq:s(ba--n)} and \eqref{eq:nn-1--1jb-1} that
\begin{align}\label{eq:S(sigma_A)j=1}
 S(\sigma_A)_j^*=& (b+1)\eta_{2,b}\eta_{b+2,n}
 -\sum_{i=4}^{b+1} S(\sigma_{\{i-1,i\},i})_1^*\sha\eta_{i+1,n} -\sum_{i=b+2}^{n-1}S(\sigma_{\{b,i\},i})_1^*\sha\eta_{i+1,n}.
\end{align}

\eqref{lem:item:j=1b=1} If $j=1$ and $b=1$, then the second summand of Eq.\,\eqref{eq:S(sigma_A)j=1} is $0$, so that  Eq.\,\eqref{eq:S(sigma_A)j=1} becomes
\begin{align*}
 S(\sigma_A)_j^*=\eta_{2,n}-\sum_{i=3}^{n-1}S(\sigma_{\{1,i\},i})_1^*\sha\eta_{i+1,n}=\eta_{2,n}+\sum_{i=3}^{n-1}(-1)^{i}(2\sha \delta_{i,4})3\sha\eta_{i+1,n},
\end{align*}
where the second equality follows from the induction hypothesis on part \eqref{lem:item:j=1b=1}.
Substituting $j=k=2$ in Lemma \ref{lem:eta2n+-1i} yields that $S(\sigma_A)_j^*=(-1)^{n+1}(2\sha \delta_{n,4})3$.

\eqref{lem:item:j=1b=2}  If $j=1$ and $b=2$, then the second summand of Eq.\,\eqref{eq:S(sigma_A)j=1} is $0$, so by the induction hypothesis on part \eqref{lem:item:j=1b=2}
\begin{align*}
 S(\sigma_A)_j^*=32\eta_{4,n}-\sum_{i=4}^{n-1}S(\sigma_{\{2,i\},i})_1^*\sha\eta_{i+1,n}
 =\delta_{3,2}\eta_{4,n}+\sum_{i=4}^{n-1}(-1)^{i+1}(32\sha\delta_{i,5})4\sha\eta_{i+1,n}.
\end{align*}
Taking $j=2$ and $k=3$ in Lemma \ref{lem:eta2n+-1i} gives that
$S(\sigma_A)_j^*=(-1)^{n} (32\sha \delta_{n,5})4$.

\eqref{lem:item:j=1bgeq3}
If $j=1$ and $3\leq b\leq n-2$, then apply Lemma \ref{lem:etaibn-1}\eqref{lem:item:sba1a-1'} and the induction hypothesis on part \eqref{lem:item:j=1bgeq3} to the second and third summands of Eq.\,\eqref{eq:S(sigma_A)j=1}, we get
\begin{align}\label{eq:j=1bgeq3pf}
 S(\sigma_A)_j^*
=&\, (b+1)\eta_{2,b}\eta_{b+2,n}-\sum_{i=4}^{b+1} i\eta_{2,i-1}\sha\eta_{i+1,n}
 +\sum_{i=4}^{b+1}\sum_{k=4}^{i-1}(-1)^{i+k+1}\left(k\eta_{2,k-1}\sha\delta_{i,k+1}\right)\sha\eta_{i+1,n}\notag\\
& \hspace{-10mm}+\sum_{i=b+2}^{n-1} (-1)^{i+b}\left[(b+1)\eta_{2,b-1}\sha\delta_{i,b+2}\right]b\sha\eta_{i+1,n} +\sum_{i=b+2}^{n-1}\sum_{k=4}^{b} (-1)^{i+k+1}\left(k\eta_{2,k-1}\sha\delta_{i,k+1}\right)\sha\eta_{i+1,n}.
\end{align}
Applying Lemma \ref{lem:etadeltaetan} with $k=b$, the fourth term on the right-hand side of Eq.\,\eqref{eq:j=1bgeq3pf} is equal to
\begin{align*}
-(b+1)\eta_{2,b}\eta_{b+2,n}+(b+1)\eta_{2,b}\sha\eta_{b+2,n}+(-1)^{n+b+1}\left[(b+1)\eta_{2,b-1}\sha\delta_{n,b+2}\right]b,
\end{align*}
while applying Lemma \ref{lem:-1ideltaij+1shaeta}, the sum of the third and fifth summands on the right-hand side of Eq.\,\eqref{eq:j=1bgeq3pf} is equal to
\begin{align*}
\sum_{k=4}^{b}\sum_{i=k+1}^{n-1} (-1)^{i+k+1}\left(k\eta_{2,k-1}\sha\delta_{i,k+1}\right)\sha\eta_{i+1,n}
 =\sum_{k=4}^{b} k\eta_{2,k-1}\sha\eta_{k+1,n}
 +\sum_{k=4}^{b} (-1)^{n+k}(k\eta_{2,k-1}\sha \delta_{n,k+1}).
\end{align*}
Substituting into Eq.\,\eqref{eq:j=1bgeq3pf} gives
\begin{align*}
 S(\sigma_A)_j^*=(-1)^{n+b+1}\left[(b+1)\eta_{2,b-1}\sha\delta_{n,b+2}\right]b+\sum_{k=4}^{b} (-1)^{n+k}(k\eta_{2,k-1}\sha \delta_{n,k+1}),
\end{align*}
so the proof of part \eqref{lem:item:j=1bgeq3} follows.

We now prove  parts \eqref{lem:item:j=2b=1a=2} and \eqref{lem:item:j=2b=1ageq3}.
Assume that $j=2$, $b=1$. Then $U_j^*=V_j^*=X_j^*=0$.
It follows from Eqs.\,\eqref{eq:ushavjeq} and \eqref{eq:s(ba--n)} that
\begin{align}\label{eq:s(ba--n)j=2b=1a=2}
 S(\sigma_A)_j^*=\,-1\sha \eta_{3,n}-\sum_{i=3}^{{\rm min}\{a,n-1\}}S(\sigma_{\{1,i\},i})_2^*\sha\eta_{i+1,n}-\sum_{i=a+1}^{n-1}S(\sigma_{\{1,a\},i})_2^*\sha\eta_{i+1,n}.
\end{align}

\eqref{lem:item:j=2b=1a=2} If $j=2$, $b=1$ and $a=2$, then the second term of Eq.\,\eqref{eq:s(ba--n)j=2b=1a=2} is $0$, and hence by the induction hypothesis on part \eqref{lem:item:j=2b=1a=2}, we have
$$S(\sigma_A)_j^*=-1\sha \eta_{3,n}-\sum_{i=3}^{n-1}S(\sigma_{\{1,2\},i})_2^*\sha\eta_{i+1,n}
=-1\sha \eta_{3,n}-\sum_{i=3}^{n-1}(-1)^i(1\sha\delta_{i,3}\sha\eta_{i+1,n})=(-1)^{n}(1\sha \delta_{n,3}),$$
where the third equality is obtained from Lemma \ref{lem:-1ideltaij+1shaeta} by putting $k=2$.

\eqref{lem:item:j=2b=1ageq3} If $j=2$, $b=1$ and $a\geq3$, then, by the induction hypothesis on part \eqref{lem:item:j=2b=1ageq3}, $S(\sigma_{\{1,a\},i})_2^*=S(\sigma_{\{1,i\},i})_2^*$ when $a+1\leq i\leq n-1$.
Hence,
\begin{align*}
 S(\sigma_A)_j^*&=-1\sha \eta_{3,n}-\sum_{i=3}^{n-1}S(\sigma_{\{1,i\},i})_2^*\sha\eta_{i+1,n}&&{\rm by\ Eq.}\,\eqref{eq:s(ba--n)j=2b=1a=2}\\
 &=-1\sha \eta_{3,n}-\sum_{i=3}^{n-1}(-1)^{i}\left(1\sha \delta_{i,3}+31\sha \delta_{i,4}\right)\sha\eta_{i+1,n}&&{\rm by\ the\ induction\ hypothesis}\\
 &=-1\sha\left(\sum_{i=2}^{n-1}(-1)^{i}\delta_{i,3}\sha\eta_{i+1,n}\right)-31\sha\left(\sum_{i=3}^{n-1}(-1)^{i}\delta_{i,4}\sha\eta_{i+1,n}\right).
\end{align*}
Putting $k\in\{2,3\}$ in Lemma \ref{lem:-1ideltaij+1shaeta} yields that $ S(\sigma_A)_j^*=(-1)^{n}\left(1\sha \delta_{n,3}+31\sha \delta_{n,4}\right)$.

\eqref{lem:item:j=2b=2} If $j=2$ and $b=2$, then $U_j^*=V_j^*=0$, and hence
\begin{align*}
 S(\sigma_A)_j^*&=-1\sha \eta_{3,n}-S(321)_2^*\sha\eta_{4,n}-\sum_{i=4}^{n-1}S(\sigma_{\{2,i\},i})_2^*\sha\eta_{i+1,n}
  &&{\rm by\ Eqs.}\,\eqref{eq:ushavjeq},\ \eqref{eq:s(ba--n)}\ {\rm and}\ \eqref{eq:nn-1--1jb-1}\\
 &=-1\sha \eta_{3,n}-\sum_{i=4}^{n-1}S(\sigma_{\{2,i\},i})_2^*\sha\eta_{i+1,n}&&{\rm by}\ S(321)_2^*=0\\
 &=-1\sha \eta_{3,n}-\sum\limits_{i=4}^{n-1}(-1)^{i+1}\left[1\sha(3\sha\delta_{i,5})4\right]\sha\eta_{i+1,n}&&{\rm by\ the\ induction\ hypothesis}\\
 &=1\sha \left(-\eta_{3,n}+\sum\limits_{i=4}^{n-1}(-1)^{i}(3\sha\delta_{i,5})4\sha\eta_{i+1,n}\right).
\end{align*}
Putting $j=k=3$ in Lemma \ref{lem:eta2n+-1i} yields that $S(\sigma_{A})_j^*=(-1)^{n+1}\left[1\sha (3\sha\delta_{n,5})4\right]$.

\eqref{lem:item:j=2bgeq3}
If $j=2$ and $b\geq 3$, then $U_j^*=V_j^*=0$. Since $S(321)_2^*=0$, it follows from Eqs.\,\eqref{eq:ushavjeq}, \eqref{eq:s(ba--n)} and \eqref{eq:nn-1--1jb-1} that
\begin{align*}
 S(\sigma_A)_j^*=&\,-1\sha \eta_{3,n}-\sum_{i=4}^{b+1} S(\sigma_{\{i-1,i\},i})_2^*\sha\eta_{i+1,n}-\sum_{i=b+2}^{n-1}S(\sigma_{\{b,i\},i})_2^*\sha\eta_{i+1,n}.
 \end{align*}
By the induction hypothesis on part \eqref{lem:item:j=2bgeq3} and Lemma \ref{lem:etaibn-1}\eqref{lem:item:sba1a-1j=2'}, we get
 \begin{align*}
 S(\sigma_A)_j^*&=-1\sha \eta_{3,n}-\sum\limits_{i=4}^{n-1}(-1)^{i+1}\left[1\sha(3\sha\delta_{i,5})4+431\sha\delta_{i,5}\right]\sha\eta_{i+1,n}\\
 &=1\sha\left(-\eta_{3,n}+\sum\limits_{i=4}^{n-1}(-1)^{i}(3\sha\delta_{i,5})4\sha\eta_{i+1,n}\right)
 +431\sha\left(\sum\limits_{i=4}^{n-1}(-1)^{i}\delta_{i,5}\sha\eta_{i+1,n}\right).
\end{align*}
So the proof follows from Lemma \ref{lem:eta2n+-1i} with $k=j=3$ and Lemma \ref{lem:-1ideltaij+1shaeta} with $k=4$.

\eqref{lem:item:3leqjleqb-1}
Assume that $3\leq j\leq b-1$. Since $b\leq n-2$, we have $3\leq j\leq n-3$, which yields that $U_j^*=V_j^*=W_j^*=0$. Applying Eqs.\,\eqref{eq:ushavjeq} and \eqref{eq:nn-1--1jb-1} to Eq.\,\eqref{eq:s(ba--n)} gives that
\begin{align*}%\label{eq:s(ba--n)3leqjleqb-1}
 S(\sigma_A)_j^*=-\sum_{i=j}^{b+1} S(\sigma_{\{i-1,i\},i})_j^*\sha\eta_{i+1,n}-\sum_{i=b+2}^{n-1} S(\sigma_{\{b,i\},i})_j^*\sha\eta_{i+1,n}.
\end{align*}
Since for $i$ with $b+2\leq i\leq n-1$, we have $3\leq j\leq i-3$, it follows from Lemma \ref{lem:etaibn-1}\eqref{lem:item:sbaa'} and
the induction hypothesis on part \eqref{lem:item:3leqjleqb-1} that $S(\sigma_{\{b,i\},i})_j^*=S(\sigma_{\{i-1,i\},i})_j^*$ for $i$ with $b+2\leq i\leq n-1$, and hence
\begin{align*}%\label{eq:s(ba--n)3leqjleqb-1}
 S(\sigma_A)_j^*=-\sum_{i=j}^{n-1} S(\sigma_{\{i-1,i\},i})_j^*\sha\eta_{i+1,n}.
\end{align*}
Thus, part \eqref{lem:item:3leqjleqb-1} is true by Lemma \ref{lem:etaibn-1}\eqref{lem:item:sbaa'} and the proof of Eq.\,\eqref{eq:s(ba--n)3leqjleqb-111}.

\eqref{lem:item:3leqj=b} Suppose $3\leq j=b$.
Then $j\leq n-2$, and we see that  $U_j^*=V_j^*=W_j^*=0$. It follows from Eqs.\,\eqref{eq:ushavjeq}, \eqref{eq:s(ba--n)} and \eqref{eq:nn-1--1jb-1} that
\begin{align*}
 S(\sigma_A)_j^*=&-\sum_{i=b}^{b+1} S(\sigma_{\{i-1,i\},i})_b^*\sha\eta_{i+1,n} -\sum_{i=b+2}^{n-1}S(\sigma_{\{b,i\},i})_b^*\sha\eta_{i+1,n}\\
 =&\, -S(\sigma_{\{j-1,j\},j})_j^*\sha\eta_{j+1,n} -\sum_{i=j+1}^{n-1}S(\sigma_{\{j,i\},i})_j^*\sha\eta_{i+1,n}.
\end{align*}
By Lemma \ref{lem:etaibn-1}\eqref{lem:item:sbaaj=b+1geq3'}, we have
$S(\sigma_{\{j-1,j\},j})_j^*\sha\eta_{j+1,n}=-\left(1\sha (j-1)\eta_{2,j-2}-\delta_{j-1,j-2}\eta_{1,j-3}\right)\sha \eta_{j+1,n}$,
and by the induction hypothesis on part \eqref{lem:item:3leqj=b}, we have
\begin{align*}
\sum_{i=j+1}^{n-1}S(\sigma_{\{j,i\},i})_j^*\sha\eta_{i+1,n}=
&\,(-1)^{j+1}\left[1\sha (j-1)\eta_{2,j-2}-\delta_{j-1,j-2}\eta_{1,j-3}\right]\sha \left(\sum_{i=j+1}^{n-1}(-1)^i\delta_{i,j+1}\sha\eta_{i+1,n}\right)\\
&\,+(-1)^{j+1}\left(1\sha (j+1)\eta_{2,j-1}\right)\sha \left(\sum_{i=j+1}^{n-1}(-1)^i\delta_{i,j+2}\sha\eta_{i+1,n}\right).
\end{align*}
Applying Lemma \ref{lem:-1ideltaij+1shaeta} with $k=j$ to the first term, while applying Lemma \ref{lem:-1ideltaij+1shaeta} with $k=j+1$ to the second term in the above equation, we obtain that
\begin{align*}
\sum_{i=j+1}^{n-1}S(\sigma_{\{j,i\},i})_j^*\sha\eta_{i+1,n}=&\,\left(1\sha (j-1)\eta_{2,j-2}-\delta_{j-1,j-2}\eta_{1,j-3}\right)\sha \left(\eta_{j+1,n}+(-1)^{n+j}\delta_{n,j+1}\right)\\
&\,+(-1)^{n+j}\left(1\sha (j+1)\eta_{2,j-1}\right)\sha \delta_{n,j+2},
\end{align*}
from which the proof follows.

\eqref{lem:item:3leqj=b+1a=b+1}
Let $3\leq j=b+1$ and $a=b+1$. Then  $U_j^*=V_j^*=W_j^*=Y_j^*=0$.
From Eqs.\,\eqref{eq:ushavjeq} and \eqref{eq:s(ba--n)} we see that
\begin{align*}%\label{eq:3leq j=b+1pf1}
S(\sigma_A)_{j}^*=-S(\sigma_{\{b,b+1\},b+1})_{b+1}^*\sha\eta_{b+2,n}
 -\sum_{i=a+1}^{n-1}S(\sigma_{\{b,a\},i})_{b+1}^*\sha\eta_{i+1,n}
 = -\sum_{i=j}^{n-1}S(\sigma_{\{j-1,j\},i})_{j}^*\sha\eta_{i+1,n}.
\end{align*}
By the induction hypothesis on part \eqref{lem:item:3leqj=b+1a=b+1}, we have
\begin{align*}
S(\sigma_A)_{j}^*&=-\sum_{i=j}^{n-1}(-1)^{i+j+1}\left[\left(1\sha (j-1)\eta_{2,j-2}-\delta_{j-1,j-2}\eta_{1,j-3}\right)\sha \delta_{i,j+1}\right]\sha\eta_{i+1,n}\\
&=(-1)^{j}\left(1\sha (j-1)\eta_{2,j-2}-\delta_{j-1,j-2}\eta_{1,j-3}\right)\sha\left(\sum_{i=j}^{n-1}(-1)^{i} \delta_{i,j+1}\sha\eta_{i+1,n}\right),
\end{align*}
and the proof follows from Lemma \ref{lem:-1ideltaij+1shaeta} with $k=j$.

\eqref{lem:item:3leqj=b+1ageqb+1}
Let $3\leq j=b+1$ and $a\geq b+2$. Then $V_j^*=W_j^*=0$. By Eqs.\,\eqref{eq:ushavjeq} and \eqref{eq:s(ba--n)}, $S(\sigma_A)_{j}^*=P+Q$, where
\begin{align*}%\label{eq:3leq j=b+1pf1}
P=&-(ab\eta_{1,b-1}\eta_{b+1,a-1}\eta_{a+1,n})_{b+1}^*,\\
Q=&-S(\sigma_{\{b,b+1\},b+1})_{b+1}^*\sha\eta_{b+2,n}-\sum_{i=b+2}^{{\rm min}\{a,n-1\}}S(\sigma_{\{b,i\},i})_{b+1}^*\sha\eta_{i+1,n}
 -\sum_{i=a+1}^{n-1}S(\sigma_{\{b,a\},i})_{b+1}^*\sha\eta_{i+1,n}.
\end{align*}
By definitions,
\begin{align*}%\label{eq:3leq j=b+1pf1}
P=\begin{cases}
-n(n-2)\eta_{1,n-3} &\text{if $a=n$ and $b= n-2$},\\
0 &\text{otherwise}.
\end{cases}
\end{align*}

First consider the case $a=n$.
If $b=n-2$, then, by Lemma \ref{lem:etaibn-1}\eqref{lem:item:sbaaj=b+1geq3'},
\begin{align*}%\label{eq:3leq j=b+1pf1}
Q=-S(\sigma_{\{n-2,n-1\},n-1})_{n-1}^*\sha n=\left(1\sha(n-2)\eta_{2,n-3}-\delta_{n-2,n-3}\eta_{1,n-4}\right)\sha n.
\end{align*}
Thus,
$$ S(\sigma_A)_{j}^*=P+Q=\left(1\sha(n-2)\eta_{2,n-3}-\delta_{n-2,n-3}\eta_{1,n-4}\right)\sha n-n(n-2)\eta_{1,n-3},
$$
agreeing with the formula given by part \eqref{lem:item:3leqj=b+1ageqb+1}.
If $b\neq n-2$, then $b\leq n-3$. Since $a=n$, there follows
\begin{align*}%\label{eq:3leq j=b+1pf}
S(\sigma_A)_{j}^*=Q=-\sum_{i=b+1}^{n-1}S(\sigma_{\{b,i\},i})_{b+1}^*\sha\eta_{i+1,n}= -\sum_{i=j}^{n-1}S(\sigma_{\{j-1,i\},i})_{j}^*\sha\eta_{i+1,n}.
\end{align*}
By Lemma \ref{lem:etaibn-1}\eqref{lem:item:sbaaj=b+1geq3'}, we have
$S(\sigma_{\{j-1,j\},j})_{j}^*=\delta_{j-1,j-2}\eta_{1,j-3}-1\sha (j-1)\eta_{2,j-2}$, and by the induction hypothesis on part \eqref{lem:item:3leqj=b+1ageqb+1},
we see that $\sum_{i=j+1}^{n-1}S(\sigma_{\{j-1,i\},i})_{j}^*\sha\eta_{i+1,n}$ is equal to
\begin{align*}%\label{eq:3leq j=b+1pf}
\sum_{i=j+1}^{n-1}(-1)^{i+j+1}\left[\left(1\sha (j-1)\eta_{2,j-2}-\delta_{j-1,j-2}\eta_{1,j-3}\right)\sha \delta_{i,j+1}-(j+1)(j-1)\eta_{1,j-2}\sha \delta_{i,j+2}\right]\sha\eta_{i+1,n}.
\end{align*}
Hence,
\begin{align*}%\label{eq:3leq j=b+1pf}
S(\sigma_A)_{j}^*
=&\,(-1)^{j}\left(1\sha (j-1)\eta_{2,j-2}-\delta_{j-1,j-2}\eta_{1,j-3}\right)\sha\left(\sum_{i=j}^{n-1}(-1)^{i} \delta_{i,j+1}\sha\eta_{i+1,n}\right)\\
&\,+(-1)^{j+1}(j+1)(j-1)\eta_{1,j-2}\sha \left( \sum_{i=j+1}^{n-1}(-1)^{i} \delta_{i,j+2}\sha\eta_{i+1,n}\right).
\end{align*}
Applying Lemma \ref{lem:-1ideltaij+1shaeta} with $k=j$ to the first term, while applying Lemma \ref{lem:-1ideltaij+1shaeta} with $k=j+1$ to the second term in the above equation,
gives the desired identity.

We now consider the case $a\leq n-1$.
Then $S(\sigma_A)_{j}^*=Q$. Since $a\geq b+2$, by the induction hypothesis on part \eqref{lem:item:3leqj=b+1ageqb+1} we have $S(\sigma_{\{b,a\},i})_{b+1}^*=S(\sigma_{\{b,i\},i})_{b+1}^*$ for all $i$ with $a+1\leq i\leq n-1$, and hence
\begin{align*}
 S(\sigma_A)_{j}^*
 =&-S(\sigma_{\{b,b+1\},b+1})_{b+1}^*\sha\eta_{b+2,n}-\sum_{i=b+2}^{a}S(\sigma_{\{b,i\},i})_{b+1}^*\sha\eta_{i+1,n}
 -\sum_{i=a+1}^{n-1}S(\sigma_{\{b,a\},i})_{b+1}^*\sha\eta_{i+1,n}\\
 =&\,-S(\sigma_{\{j-1,j\},j})_{j}^*\sha\eta_{j+1,n} -\sum_{i=j+1}^{n-1}S(\sigma_{\{j-1,a\},i})_{j}^*\sha\eta_{i+1,n}.
\end{align*}
By Lemma \ref{lem:etaibn-1}\eqref{lem:item:sbaaj=b+1geq3'} and the induction hypothesis on part \eqref{lem:item:3leqj=b+1ageqb+1},
\begin{align*}
 S(\sigma_A)_{j}^*=&\,\left(1\sha (j-1)\eta_{2,j-2}-\delta_{j-1,j-2}\eta_{1,j-3}\right)\sha\eta_{j+1,n} \\
 &+(-1)^{j}\left(1\sha (j-1)\eta_{2,j-2}-\delta_{j-1,j-2}\eta_{1,j-3}\right)\sha \left(\sum_{i=j+1}^{n-1}(-1)^{i}\delta_{i,j+1}\sha\eta_{i+1,n}\right)\\
&\hspace{30mm}+(-1)^{j+1}(j+1)(j-1)\eta_{1,j-2}\sha \left(\sum_{i=j+1}^{n-1}(-1)^{i}\delta_{i,j+2}\sha\eta_{i+1,n}\right).
\end{align*}
Applying Lemma \ref{lem:-1ideltaij+1shaeta} with $k=j$ to the sum of the first and second terms, while applying Lemma \ref{lem:-1ideltaij+1shaeta} with $k=j+1$ to the third term in the above equation, we obtain the desired formula given by  part \eqref{lem:item:3leqj=b+1ageqb+1}.

\eqref{lem:item:b+2leqjleqnjleqa-1} Suppose $b+2\leq j\leq n$ and $j\leq a-1$. Then $3\leq j\leq n-1$ and $a\geq b+3$, so $V_j^*=W_j^*=X_j^*=0$. It follows from Eqs.\,\eqref{eq:ushavjeq} and \eqref{eq:s(ba--n)} and the induction hypothesis on part \eqref{lem:item:b+2leqjleqnjleqa-1} that
\begin{align*}
S(\sigma_A)_{j}^*&=-(ab\eta_{1,b-1}\eta_{b+1,a-1}\eta_{a+1,n})_j^*-\sum_{i=j}^{{\rm min}\{a,n-1\}}S(\sigma_{\{b,i\},i})_{j}^*\sha\eta_{i+1,n}-\sum_{i=a+1}^{n-1}S(\sigma_{\{b,a\},i})_{j}^*\sha\eta_{i+1,n}\\
&=-(ab\eta_{1,b-1}\eta_{b+1,a-1}\eta_{a+1,n})_j^*-S(\sigma_{\{b,j\},j})_{j}^*\sha\eta_{j+1,n}-\sum_{i=j+1}^{n-1}S(\sigma_{\{b,a\},i})_{j}^*\sha\eta_{i+1,n}.
\end{align*}
If $a=n$ and $j=n-1$, then $b\leq j-2=n-3$. By the induction hypothesis on part \eqref{lem:item:b+2leqjleqnj=a}, we see that
\begin{align*}
S(\sigma_A)_{j}^*
&=-nb\eta_{1,b-1}\eta_{b+1,n-2}-S(\sigma_{\{b,n-1\},n-1})_{n-1}^*\sha n\\
&=-nb\eta_{1,b-1}\eta_{b+1,n-2}-\left[(n-2)b \eta_{1,b-1}\eta_{b+1,n-3}\sha \delta_{n-1,n}\right]\sha n\\
&=-(n-2)b \eta_{1,b-1}\eta_{b+1,n-3}\sha n-nb\eta_{1,b-1}\eta_{b+1,n-2},
\end{align*}
which coincides with the identity given by part \eqref{lem:item:b+2leqjleqnjleqa-1}.
If either $a\neq n$ or $j\neq n-1$, then we have $(ab\eta_{1,b-1}\eta_{b+1,a-1}\eta_{a+1,n})_j^*=0$. From the induction hypothesis on part \eqref{lem:item:b+2leqjleqnj=a} we conclude that $S(\sigma_{\{b,j\},j})_{j}^*=(j-1)b \eta_{1,b-1}\eta_{b+1,j-2}$,
which together with the induction hypothesis on part \eqref{lem:item:b+2leqjleqnjleqa-1} yields that
\begin{align*}
S(\sigma_A)_{j}^*
=&\,-(j-1)b \eta_{1,b-1}\eta_{b+1,j-2}\sha\eta_{j+1,n}\\
&\,-\sum_{i=j+1}^{n-1}(-1)^{i+j}\left[(j-1)b \eta_{1,b-1}\eta_{b+1,j-2}\sha \delta_{i,j+1}+(j+1)b\eta_{1,b-1}\eta_{b+1,j-1}\sha \delta_{i,j+2}\right]\sha\eta_{i+1,n}\\
=&\,(-1)^{j+1}(j-1)b \eta_{1,b-1}\eta_{b+1,j-2}\sha\left(\sum_{i=j}^{n-1}(-1)^{i}\delta_{i,j+1}\sha\eta_{i+1,n}\right)\\
&\,\hspace{30mm}+(-1)^{j+1}(j+1)b\eta_{1,b-1}\eta_{b+1,j-1}\sha\left(\sum_{i=j+1}^{n-1}(-1)^{i}\delta_{i,j+2}\sha\eta_{i+1,n}\right).
\end{align*}
Using Lemma \ref{lem:-1ideltaij+1shaeta} with $k=j$ to the first term, while using Lemma \ref{lem:-1ideltaij+1shaeta} with $k=j+1$ to the second term in the above equation, we obtain the desired result.

\eqref{lem:item:b+2leqjleqnj=a} Suppose $b+2\leq j\leq n$ and $j=a$. If $j=n$, then $U_j^*=Z_j^*=0$. From Eqs.\,\eqref{eq:ushavjeq} and \eqref{eq:s(ba--n)} we see that
\begin{align*}
S(\sigma_A)_{j}^*=&1\sha(b+1)\eta_{2,b}\eta_{b+2,n-1}-12\sha\eta_{3,n-1}-\sum_{i=3}^{b+1}S(\sigma_{\{i-1,i\},i})\sha\eta_{i+1,n-1}-\sum_{i=b+2}^{n-1}S(\sigma_{\{b,i\},i})\sha\eta_{i+1,n-1}.
\end{align*}
On the other hand, by Eq.\eqref{eq:antipodeformula1},
\begin{align*}
S(\sigma_{\{b,n-1\},n-1})=&-(n-1)b\eta_{1,b-1}\eta_{b+1,n-2}+1\sha (b+1)\eta_{2,b}\eta_{b+2,n-1}-12\sha\eta_{3,n-1}\\
&-\sum_{i=3}^{b+1}S(\sigma_{\{i-1,i\},i})\sha\eta_{i+1,n-1}-\sum_{i=b+2}^{n-2}S(\sigma_{\{b,i\},i})\sha\eta_{i+1,n-1},
\end{align*}
so that
\begin{align*}
S(\sigma_A)_{j}^*=(n-1)b\eta_{1,b-1}\eta_{b+1,n-2}=(-1)^{n+j}\left[(j-1)b \eta_{1,b-1}\eta_{b+1,j-2}\sha \delta_{n,j+1}\right].
\end{align*}

If $j=n-1$, then $U_j^*=V_j^*=W_j^*=X_j^*=Z_j^*=0$. By Eqs.\,\eqref{eq:ushavjeq} and \eqref{eq:s(ba--n)}, together with the induction hypothesis on part \eqref{lem:item:b+2leqjleqnj=a}, we have
\begin{align*}
S(\sigma_A)_{j}^*=&-S(\sigma_{\{b,n-1\},n-1})_{n-1}^*\sha n=-S(\sigma_{\{b,j\},j})_{j}^*\sha n=(-1)^{n+j}\left[(j-1)b \eta_{1,b-1}\eta_{b+1,j-2}\sha \delta_{n,j+1}\right].
\end{align*}

If $b+2\leq j\leq n-2$, then since $j=a$, we have $U_j^*=V_j^*=W_j^*=X_j^*=0$. By Eqs.\,\eqref{eq:ushavjeq} and \eqref{eq:s(ba--n)}, we see that
\begin{align*}
S(\sigma_A)_{j}^*=-S(\sigma_{\{b,a\},a})_{j}^*\sha\eta_{a+1,n}-\sum_{i=a+1}^{n-1}S(\sigma_{\{b,a\},i})_{j}^*\sha\eta_{i+1,n}
=-\sum_{i=j}^{n-1}S(\sigma_{\{b,a\},i})_{j}^*\sha\eta_{i+1,n}.
\end{align*}
By the induction hypothesis on part \eqref{lem:item:b+2leqjleqnj=a} and Lemma \ref{lem:-1ideltaij+1shaeta} with $k=j$, we obtain that
\begin{align*}
S(\sigma_A)_{j}^*=&(-1)^{j+1}(j-1)b \eta_{1,b-1}\eta_{b+1,j-2}\sha \left(\sum_{i=j}^{n-1}(-1)^{i}\delta_{i,j+1} \sha \eta_{i+1,n}\right)\\
=&(-1)^{n+j}\left[(j-1)b\eta_{1,b-1}\eta_{b+1,j-2}\sha \delta_{n,j+1}\right].
\end{align*}

\eqref{lem:item:b+2leqjleqnjgeqa+1} This follows from Lemma \ref{lem:etai-n}.
\end{proof}

The Conjectures \ref{conj1} and \ref{conj2} now follow directly from Theorem \ref{maintheorem}.

\begin{coro}
Conjectures \ref{conj1} and \ref{conj2} are true.
\end{coro}
\begin{proof}
For $A=\{a\}$ where $1<a\leq n$, we have $\sigma_{A}=\sigma_{\{1,a\}}$. According to Lemma \ref{lem:etai-n}, we have
\begin{align}\label{eq:s(siga1)}
     S(\sigma_{A})= \sum_{j=1}^{a}S(\sigma_{A})_j.
\end{align}
Putting $b=1$, it follows from Theorem \ref{maintheorem}\eqref{lem:item:j=1b=1} that $S(\sigma_{A})_1 = (-1)^{n-1}(2 \sha \delta_{n,4})31$.
By Theorem \ref{maintheorem}\eqref{lem:item:j=2b=1a=2} and \eqref{lem:item:b+2leqjleqnj=a}, we have
\begin{align*}
 S(\sigma_{A})_a&=\begin{cases}
 (-1)^{n}\left(1 \sha \delta_{n,3}\right)a & {\rm if}\ a=2,\\
 (-1)^{n+a}\left[(a-1)\eta_{1,a-2} \sha \delta_{n,a+1}\right]a &{\rm if}\ a\geq 3
 \end{cases}\\
 &=(-1)^{n+a}\left[(a-1)\eta_{1,a-2} \sha \delta_{n,a+1}\right]a.
\end{align*}
For $2\leq j\leq a-1$, it follows from Theorem \ref{maintheorem}\eqref{lem:item:j=2b=1ageq3} and \eqref{lem:item:b+2leqjleqnjleqa-1} that
\begin{align*}
 S(\sigma_{A})_j &=\begin{cases}
 (-1)^{n}\left(1\sha \delta_{n,3}+31\sha \delta_{n,4}\right)j& {\rm if}\ j=2,\\
 (-1)^{n+j}[(j-1)\eta_{1,j-2} \sha \delta_{n,j+1} + (j+1)\eta_{1,j-1} \sha \delta_{n,j+2}]j& {\rm if}\ 3\leq j\leq a-1
 \end{cases}\\
 &=(-1)^{n+j}[((j-1)\eta_{1,j-2} \sha \delta_{n,j+1})j + ((j+1)\eta_{1,j-1} \sha \delta_{n,j+2})j],
\end{align*}
so the proof of Conjecture \ref{conj1} follows.

Now let $A=\{a, 2\}$ with $2 < a \leq n$. Putting $b=2$ in Theorem \ref{maintheorem}, then Eq.~\eqref{eq:s(siga1)} still holds, but now it follows from
Theorem \ref{maintheorem}\eqref{lem:item:j=1b=2} and \eqref{lem:item:j=2b=2} that
\begin{align*}
 S(\sigma_{A})_1 =& (-1)^{n}(32 \sha \delta_{n,5})41,\quad
 S(\sigma_{A})_2 =(-1)^{n-1}\left[1\sha(3\sha\delta_{n,5})4\right]2.
\end{align*}

If $a=3$, then, by Theorem \ref{maintheorem}\eqref{lem:item:3leqj=b+1a=b+1},
\begin{align*}
 S(\sigma_{A})_3=(-1)^{n}\left[\left(1\sha 2\eta_{2,1}-\delta_{2,1}\eta_{1,0}\right)\sha \delta_{n,4}\right]3
 =(-1)^{n}(12\sha \delta_{n,4})3.
\end{align*}
Hence in this case Conjecture \ref{conj2} is true.

If $a\geq4$, then it follows from Theorem \ref{maintheorem}\eqref{lem:item:3leqj=b+1ageqb+1} that
\begin{align*}
 S(\sigma_{A})_3 =&(-1)^{n}\left[\left(1\sha 2\eta_{2,1}-\delta_{2,1}\eta_{1,0}\right)\sha \delta_{n,4}-42\eta_{1,1}\sha \delta_{n,5}\right]3
 =(-1)^{n}\left(12\sha \delta_{n,4}-421\sha \delta_{n,5}\right)3.
\end{align*}
According to Theorem \ref{maintheorem}\eqref{lem:item:b+2leqjleqnjleqa-1} and \eqref{lem:item:b+2leqjleqnj=a}, for any $j\geq 4$, we have
\begin{align*}
    S(\sigma_{A})_j =\begin{cases}
    (-1)^{n+j} [(j-1)21\eta_{3,j-2} \sha \delta_{n,j+1} ]j&{\rm if}\ j=a,\\
    (-1)^{n+j}[(j+1)21\eta_{3,j-1} \sha \delta_{n,j+2}]j + (-1)^{n+j}[(j-1)21\eta_{3,j-2} \sha \delta_{n,j+1}]j&{\rm if}\ 4\leq j\leq a-1.
    \end{cases}
\end{align*}
Summing over $j$ gives that
\begin{align*}
 S(\sigma_{A}) =& (-1)^{n}\left[\left(32 \sha \delta_{n,5}\right)41 + \left(12 \sha \delta_{n,4}\right)3\right] + (-1)^{n-1}\left[\left(1 \sha (3 \sha \delta_{n,5}\right)4)2\right]\\
 &+ \sum_{j=3}^{a-1} (-1)^{n+j}\left[(j+1)21\eta_{3,j-1} \sha \delta_{n,j+2}\right]j+ \sum_{j=4}^{a}(-1)^{n+j}\left[(j-1)21\eta_{3,j-2} \sha \delta_{n,j+1}\right]j.
\end{align*}
So the proof for $a\geq 4$ follows from substituting $j$ by $j+1$ in the fourth summand. Therefore, Conjecture \ref{conj2} is true.
\end{proof}

\noindent
{\bf Acknowledgements.}
We would like to thank the anonymous referees for their very careful reading and extremely valuable comments that helped to greatly improve the exposition.
This work was partially supported by the National Natural Science Foundation of China (Grant Nos. 12071377 and 12071383).

\end{document}